\documentclass[a4paper,11pt,reqno]{amsart}
\usepackage[utf8]{inputenc}
\usepackage[T1]{fontenc}
\usepackage{lmodern}
\usepackage[english]{babel}
\usepackage{amsmath,a4wide,wasysym}
\usepackage{xfrac}
\usepackage{esint}
\usepackage{stmaryrd,mathrsfs,bm,amsthm,mathtools,yfonts,amssymb,color,braket,booktabs,graphicx,graphics,amsfonts}
\usepackage{latexsym,microtype,indentfirst,hyperref}
\usepackage{xcolor}
\usepackage{courier}

\begingroup
\newtheorem{theorem}{Theorem}[section]

\newtheorem{corollary}[theorem]{Corollary}
\endgroup

\theoremstyle{definition}
\newtheorem{definition}[theorem]{Definition}

\newtheorem{conjecture}[theorem]{Conjecture}

\newcommand\bB{\mathbf{B}}
\newcommand\Sing{\mathrm{Sing}}

\numberwithin{equation}{section}

\begin{document}

%\volumetitle{ICM 2022} % Don't alter this line.

\title[Regularity for the area functional]{The regularity theory for the area functional (in geometric measure theory)}

\author{Camillo De Lellis}
\address{School of Mathematics, Institute for Advanced Study, 1 Einstein Dr., Princeton NJ 05840, USA\\
and Universit\"at Z\"urich}
\email{camillo.delellis@math.ias.edu}

\begin{abstract}
The aim of this article is to give a rather extensive, and yet nontechnical, account of the birth of the regularity theory for generalized minimal surfaces, of its various ramifications along the decades, of the most recent developments, and of some of the remaining challenges.
\end{abstract}

\maketitle

\section{Introduction}

Let $U$ be a bounded open subset of the Euclidean space $\mathbb R^{m+n}$ and let $\Sigma\subset U$ be an $m$-dimensional surface (e.g. a $C^1$ $m$-dimensional submanifold, but we will allow more general concept of surfaces for most of this note). Then $\Sigma$ is said to be a ``critical point of the area functional'', or more commonly a ``minimal surface'', if
\begin{equation}\label{e:critical}
\left. \frac{d}{dt}\right|_{t=0} \textrm{Vol}^m (\Phi_t (\Sigma)) = 0
\end{equation}
for every smooth one-parameter family of diffeomorphisms $[-\delta, \delta] \ni t \mapsto \Phi_t$ of $\overline{U}$ such that:
\begin{itemize}
\item[(a)] $\Phi_t (x)=x$ for every $x\in \partial U$ and every $t$;
\item[(b)] $\Phi_0 (x)=x$ for every $x\in \overline{U}$. 
\end{itemize}
Here $\textrm{Vol}^m$ denotes a suitable concept of $m$-dimensional volume: in the case of classical submanifolds we can take the usual one from differential geomery.
 
A notable example of minimal surfaces are those that minimize the volume in some suitable class $\mathscr{C}$. It suffices to assume that $\mathscr{C}$ is closed under deformations satisfying (a) and (b) above to conclude that any minimizer in $\mathscr{C}$ is necessarily a critical point of the area functional. 

We call the attention of the reader to condition (a): the deformations fix the boundary of the open set $U$. Thus an example of a class $\mathscr{C}$ is that formed by those surfaces $\Sigma$'s whose boundary (in a suitable sense, for instance we can take the usual one of differential topology, if we are dealing with smooth surfaces) is a fixed one $\Gamma$ contained in $\partial U$. Such minimizer is then a surface of ``least area spanning the contour $\Gamma$''. However we ultimately have to agree on the very definition of ``admissible surface'' (Must it be embedded or do we allow self intersections? Do we allow any topological type? In fact must it be smooth or do we allow singularities? If we allow singularities, which type should we allow?), on what it means to span $\Gamma$, and how we define its volume.

Having assumed that we have answered all the above questions, i.e. that we have selected a suitable class $\mathscr{C}$ and a related concept of $m$-dimensional volume, a minimizer in $\mathscr{C}$ can be regarded as one possible solution to a celebrated problem in the calculus of variations, which goes under the name of {\em Plateau problem}. Indeed the Belgian physicist Joseph Plateau investigated it in the early XIX century with the intention of finding a good description of soap films. However the problem had already appeared in the mathematical literature decades before the investigations of Plateau, and can be found in the works of Lagrange, Meusnier, Monge and L\'egendre. In particular Lagrange considered minimizers of the area as early as the 1760s and he used his newly established method (which leads to what nowadays are called ``Euler-Lagrange'' conditions for minima of integral energies (cf. \cite{Lagrange})) to describe $2$-dimensional minimal graphs in $\mathbb R^3$ through a suitable partial differential equation. 

As it is well-known, if $\Sigma$ is of class $C^2$, minimality in the sense of \eqref{e:critical} is equivalent to the vanishing of the mean curvature vector. The latter is a condition which can be explained without any knowledge of differential geometry and it is in fact fairly easy to describe to anybody with a basic knowledge of multivariable calculus. Having fixed a point $p_0\in \Sigma$, choose first an orthonormal system of coordinates so that $\Sigma$ is the graph of a map $\psi: \mathbb R^m \supset \Omega \to \mathbb R^n$ with the properties that 
\begin{itemize}
\item $p_0 = (0, \psi (0))$ (i.e. $p_0$ is the origin of the coordinate system);
\item and $\nabla \psi (0)=0$ (i.e. the tangent to $\Sigma$ at the origin is horizontal).
\end{itemize} 
Then the mean curvature vector of $\Sigma$ vanishes at $p_0$ if and only if
 $\Delta \psi (0) =0$. One way to think about minimal surfaces is thus to understand them as solutions to a (somewhat complicated) nonlinear elliptic system of partial differential equations which linearize to the Laplace equation, namely $\Delta \psi = 0$, when we rotate the coordinates so that the tangent to the graph is horizontal.
 
The Laplace equation is universally considered as the prototypical elliptic partial differential equation of second order and its solutions, i.e. the harmonic functions, are the prototype to understand the behavior of solutions to more general second order elliptic PDEs. The Laplace equation as well has a variational flavor, since it characterizes critical points of the Dirichlet energy $\int |\nabla \psi|^2$. But one could argue that minimal surfaces are even more natural objects than harmonic functions: indeed a surface is minimal independently of the system of coordinates used to describe the ambient Euclidean space, while rigid motions of graphs which ``mix domain and target'' do not preserve the harmonicity: the latter is a concept which depends strongly on the selection of the dependent and independent variables used to describe the surface as a graph.

Critical points of the area functional have fascinated (and have been the object of study of) generations of mathematicians throughout at least two centuries and a half. One very interesting aspect of minimal surface theory is that it is relatively easy to produce singular examples. A particularly simple instance is given by holomorphic subvarieties in $\mathbb C^n$: if we identify $\mathbb C^n$ with $\mathbb R^{2n}$ and we understand holomorphic subvarieties of dimension $k$ as $2k$-dimensional surfaces (with singularities), then the latter are always minimal. In fact they are much more than just minimal: they minimize the area among a vast class of possible deformations. So an object as seemingly innocent as the complex algebraic curve $\Sigma = \{(z,w)\in \mathbb C^2: z^2=w^3\}$ is a $2$-dimensional minimal surface in $\mathbb R^4$ and in fact it is minimal according to one of the most restrictive meanings that we can give. However the origin is a point where $\Sigma$ is not a regular submanifold: in particular there is no neighborhood of $0$ in which it can be described as the graph of a function (at least if we understand our functions ``classically'' and do not allow them to take more than one value at each fixed point of their domain). 

Another simple example is given by the connected set $E$ of least length which contains three noncolinear points $p_1, p_2, p_3\in \mathbb R^2$. Such set is the union of three distinct segments $\sigma_i$'s which:
\begin{itemize}
\item have $p_i$ as one endpoint,
\item have all a common $q$ as the other endpoint, 
\item meet at $q$ forming angles of 120$^\circ$ degrees.
\end{itemize}
Again $q$ is a ``singular point'' in the sense of differential topology: $E$ is not a submanifold of $\mathbb R^2$ in any neighborhood of $q$.

 If we accept that a good theory of minimal surfaces must include singular objects, we then open a Pandora box, in particular the following seemingly innocent questions immediately come to mind:
\begin{itemize}
\item How do we define the admissible surfaces, or otherwise put, which kind of singularities do we allow?
\item Which kind of deformations do we take into account?
\item What is the $m$-dimensional volume of a singular surface?
\end{itemize}
All these questions can be studied from different points of view and can be given very different answers depending on which goals one has in mind. For instance answers which deal efficiently with the problem of minimizing the $m$-dimensional volume of surfaces in some fixed homology class of a given Riemannian manifold do not seem to give a satisfactory description of the complexity of soap films in real life. On the other hand, even though real life soap films display singularities, it can be proved that any $2$-dimensional integral homology class in a closed smooth Riemannian $3$-manifold has a smooth representative which minimizes the area. At any rate, whichever the goal, a rather large number of answers to these questions can be given in a subject of modern mathematics called Geometric measure theory.

Geometric measure theory provides powerful tools to study various variational questions linked to the theory of minimal surfaces and has produced, in more than half a century, several notions of ``singular minimal surfaces''.  In what follows I will address several of them and, for lack of a better term, I will call all of them ``generalized minimal surfaces''.  A subtopic of Geometric measure theory, which is commonly called ``regularity theory'', studies natural questions like
\begin{itemize}
\item Under which conditions singularities can be ruled out, i.e. the generalized minimal surfaces of a particular class end up being classical minimal submanifolds?
\item How large can the set of singularities be when its existence cannot be completely ruled out?
\item Which structural properties can the singularities have? 
\end{itemize} 
The ambition of this article is to give a rather extensive, and yet nontechnical, account of the birth of this topic, of its various ramifications along the decades, of the most recent developments, and of some of the remaining challenges. Since the topic is vast and complicated, I will probably not do a good service to much of the existing literature and I emphasize from the start that I consider my views quite biased. 

\subsection{Acknowledgments}
I am very grateful to Guido De Philippis and Luca Spolaor for carefully reading a very preliminary version of this manuscript, suggesting several precious improvements, and reminding me of a few pertinent results in the literature. 
This work was partially supported by the National Science Foundation through  the  grant  FRG-1854147.

\section{Plateau's problem, criticality, and stability}

Before coming to a description of the ``regularity theory'' I will first introduce, in this section, some of the most common notions of ``generalized minimal surfaces'' considered in geometric measure theory. 

\subsection{Plateau's problem: two general approaches} As already mentioned, the Plateau problem can be loosely described as ``looking for the surfaces $\Sigma$ of least volume spanning a given contour $\Gamma$''. If the surfaces $\Sigma$ in question are $C^1$ submanifolds, it is commonly understood that the $m$-dimensional volume is the usual one from calculus books. If the given contour is as well a submanifold, a seemingly natural possibility to give a rigorous definition of ``$\Sigma$ spans $\Gamma$'' is to say that $\Gamma$ is the boundary of $\Sigma$ in the usual sense of differential topology. From this point of view it is also natural to consider ambient spaces more general than the Euclidean one, and a common natural choice is to have a general complete and smooth Riemannian ambient manifold. Practically all the ``positive results'' which we will discuss in this note have a generalization to smooth Riemannian ambient manifolds, but in order to be as nontechnical as possible, I will refrain to state the general theorems and always assume that the ambient is Euclidean. There are however some counterexamples which have been thus far found only in ambient Riemannian manifolds, and given their relevance in some of the problems examined below, I felt that they should be discussed. 

If we want to enlarge the class of surfaces (and in particular allow minima with singularities in our class $\mathscr{C}$) we then have to specify at the same time what we mean by a surface, its volume, and the fact that it spans a given contour $\Gamma$. First of all we will fix the convention that the dimension of the surfaces in $\mathscr{C}$ is $m$, the dimension of $\Gamma$ is $m-1$, and the ambient Euclidean space (or Riemannian manifold, when the general case will be discussed) has dimension $m+n$. Secondly we will restrict our attention to regular contours $\Gamma$: even though this can be relaxed considerably and one can fix certain type of nonsmooth boundaries (depending on the framework), clearly when dealing with ``boundary regularity'' theorems it is natural to assume that $\Gamma$ itself has some regularity to start with. Once we have defined our generalized class of surfaces, their generalized volume, and what it means for them to span $\Gamma$, we will say that we have a ``variational framework'' for the Plateau problem.

It is possible to subdivide the various variational frameworks proposed in geometric measure theory in two large classes, which follow two rather different philosophical approaches. I will loosely describe them as
\begin{itemize}
\item \textbf{Set-theoretic}. We insist in this case that our generalized surfaces $\Sigma$ are just merely closed subsets of the ambient space which include $\Gamma$ as a subset. The fact that they ``span'' $\Gamma$ will then be encoded in some topological condition which  $\Sigma$ must satisfy, while the $m$-dimensional Volume is defined by a suitable ``measure'' which satisfies the usual requirements of measure theory and coincides with the classical $m$-dimensional volume when $\Sigma$ is a subset of a $C^1$ surface (or a countable union of subsets of $C^1$ surfaces). 
\item \textbf{Functional-analytic}. In this case we focus first on some nice, sufficiently regular, class of surfaces $\Sigma$ and we prescribe that their boundary is $\Gamma$ in some suitable convenient sense coming from algebraic topology: let us denote this privileged class  by $\mathscr{R}$. On $\mathscr{R}$ the concept of volume will be also given in terms of classical differential geometry and algebraic topology, and this will give us a functional $\mathcal{A}$ (the area functional) on $\mathscr{R}$. We then introduce some topology on $\mathscr{R}$, for instance a distance, and the class $\mathscr{C}$ will be a suitable completion of this topological space, while the functional $\mathcal{A}$ will be extended to $\mathscr{C}$ in some natural way (for instance we can take its lower semicontinuous envelope). 
\end{itemize}
Observe that the ``classical'' parametric approach of Douglas and Rado does not fit in any of these two broad descriptions. The fact that I am not including it in the scope of these notes does not reflect any judgment on its mathematical interest: the ``classical parametric theory'' is a beautiful piece of mathematics, but it has a rather different flavor compared to the results and problems which will be discussed here.

In all the variational frameworks which we will examine, no matter whether they fall in one class or the other, there is a common and recurrent use of two important objects from geometric measure theory: the Hausdorff $m$-dimensional measure and rectifiable sets. The Hausdorff $m$-dimensional measure, which we will denote by $\mathcal{H}^m$, is a very natural way of extending the classical notion of $m$-dimensional volume to {\em any} subset of the Euclidean space (or more generally of a metric space). In fact it is just one possibility, while a general theory of such extensions can be given in terms of the so-called ``Caratheodory construction'' and we refer the reader to some of the several textbook in the literature which treats it (cf. \cite{EG,Federer,Mattila}). Rectifiable $m$-dimensional sets are a very natural class of sets 
which contain $C^1$ surfaces but is closed under many more operations which are natural from the point of view of measure theory: they consist of countable unions of closed subsets of $C^1$  $m$-dimensional submanifold plus a set of zero $\mathcal{H}^m$ measure (some people, for instance Misha Gromov, consider the latter a somewhat very unpleasant technical addition, and the author agrees that there might be an efficient theory which works without the annoying technicality of adding null sets; however most people in analysis grow accustomed to it, as ``completing'' a $\sigma$-algebra by adding sets of measure zero is a fairly common operation). 

It was a major discovery of Besicovitch in the first half of the XXth century that any set of finite $\mathcal{H}^1$ measure can be decomposed into the union of a ``rectifiable portion'' and a ``purely unrectifiable portion'' (cf. \cite{Bes1,Bes2,Bes3}). The latter is somewhat ``orthogonal'' to any $C^1$ submanifold, in the sense that intersects any $C^1$ curve in a set of $\mathcal{H}^1$ measure zero, even though it might have positive $\mathcal{H}^1$ measure. The more general theory of $m$-dimensional rectifiable sets was developed later Federer, cf. \cite{Federer1947} (and see also \cite{MR}).
 While the rectifiable part has much of the features of $C^1$ submanifolds, and can be considered as a weak version of the latter, the purely unrectifiable part behaves in a rather counterintuitive way and in what follows we will discount it: it is however one of the major early developments of geometric measure theory that one can, without loss of generality, discard unrectifiable sets pretty much in all variational theories for the area functional and I am hiding quite deep and beautiful theorems here.

\subsection{Examples of set-theoretic approaches}\label{ss:set-theoretic} The first to pioneer what I dubbed ``set-theoretic approach'' is Reifenberg in \cite{Reif1}. In his variational framework the definition of ``$E$ spans 
$\Gamma$'' is that $\Gamma$ is trivial in the relative Cech homology of $E$ (cf. \cite{Reif1} for the precise definition). More recently Harrison (cf. \cite{Harrison,HP}) suggested another, very elegant, possible definition of ``$E$ spans $\Gamma$'' which for simplicity we describe in the easiest case of $m-1$-dimensional $\Gamma$'s in $\mathbb R^{m+1}$: any closed curve $\gamma \subset \mathbb R^{m+1}\setminus \Gamma$ which is not contractible in $\mathbb R^{m+1}\setminus \Gamma$ must intersect $E$. 

Another point of view is that taken by Almgren in his theory of $(M, \varepsilon, \delta)$-minimal sets, cf. \cite{Alm76}: rather than giving a precise notion of ``spanning'' we focus on which deformations are allowed and assume that our class $\mathscr{C}$ is closed under the latter deformations. In his work Almgren gave a far-reaching existence and regularity theory, and the existence part was recently revisited and extended in \cite{FK}. Concerning deformations a very interesting point raised only recently by David is that, in practically all the works in the literature thus far, the authors used deformations which completely ``fix'' the boundary $\Gamma$, while it would be much more natural to impose that they in fact map $\Gamma$ onto itself in some controlled way (for instance they are isotopic to the identity within the class of diffeomorphisms of $\Gamma$): this idea is at the base of his recent theory of ``sliding minimizers'', cf. \cite{David1,David2,David3}. 

In all these variational frameworks for any given sequence $E_k$ of compact sets there is a natural notion of convergence, the one of Hausdorff, for which we can extract a converging subsequence. However the Hausdorff measure $\mathcal{H}^m$ does not behave well in terms of the latter convergence, in the sense that it is not lower semicontinuous. On the other hand one can suitably adjust {\em minimizing sequences} so to achieve the lower semicontinuity of $\mathcal{H}^m$: it thus suffices to prove that the limit is in the considered class $\mathscr{C}$ to achieve a minimizer. The author, in a joint work with F. Ghiraldin and F. Maggi in \cite{DGM}, pointed out that there is in fact no need to adjust the minimizing sequence and that a suitable compactness and lower-semicontinuity statement is valid for any minimizing sequence in a class $\mathscr{C}$ as soon as it allows a rather limited number of basic competitors. In particular this gives a unified framework which treats all known examples of set-theoretic approaches put forward thus far, cf. also \cite{DDG,DPDG,DPDG2}. 

From the point of view of differential and algebraic topology, all the set-theoretic approaches have some very undesirable properties. Typical set-theoretic minimizers of the Plateau problem will always have singularities: if the boundary $\Gamma$ is complicated it is energetically convenient to form ``triple junctions'' along a singularity of codimension $1$. On the other hand one of  the biggest achievements of the functional-analytic approach is that for every smooth closed embedded curve $\Gamma$ in $\mathbb R^3$ there is always a smooth oriented $2$-dimensional submanifold with boundary $\Gamma$ which minimizes the $2$-dimensional area among all smooth oriented $2$-dimensional submanifolds with boundary $\Gamma$. Likewise it is possible to show that every $2$-dimensional integral homology class in a closed Riemannian $3$-manifold has a smooth representative which minimizes the area. The set-theoretic approaches are not able to detect these two beautiful phenomena. 

On the other hand, actual soap films do form triple junctions singularities (and even more complicated ones) in real life, and these phenomena do not seem to be efficiently captured by functional-analytic frameworks (even though these type of singularities do occur in some specific situations, see below). Much of the research in the set-theoretic frameworks is thus motivated by the original intention of Plateau of finding a good variational description of soap films. In that respect the recent paper \cite{MSS} by Maggi, Scardicchio, and Stuvard pointed out that much of the investigations in the mathematical literature have thus far ignored some very relevant physical attributes of real-life soap films. Combining some of the aspects of the set-theoretic approaches with other modern techniques, like $\Gamma$-convergence, and with more accurate considerations from mathematical physics, the papers \cite{KMS,KMS2,KMS3} propose a new variational theory which promises to provide a much more accurate description of real-life soap films.

\subsection{Functional-analytic frameworks} The pioneer of functional-analytic frameworks seems to be Renato Caccioppoli. In his works \cite{C1,C2} Caccioppoli proposed the following definition of ``perimeter'' of a general (Lebesgue measurable) set of $\mathbb R^{m+1}$ (I will actually describe a slight variation of Caccioppoli's approach, but the actual differences are just of technical nature and for the purposes of this discussion I will ignore them). First of all, if the set has a $C^1$ boundary, its perimeter is defined to be the usual $m$-dimensional volume of the boundary. Next, given a general Lebesgue measurable set $E\subset \mathbb R^{m+1}$, we consider all possible sequences $E_k$ of sets with $C^1$ boundaries with the property that the Lebesgue measure of the symmetric difference $E_k \Delta E$ goes to $0$. We then consider 
\[
\liminf_{k\to\infty} \mathcal{H}^m (E_k)
\]
and we further take the infimum of all such numbers among all approximating sequences $\{E_k\}$. The latter is defined to be the perimeter of $E$. If it is finite, $E$ is commonly called {\em set of finite perimeter} or (especially if you are Italian!) {\em Caccioppoli set}. 

Caccioppoli's approach is very natural in the Calculus of Variations. We start from a class of ``good'' objects, the open sets with smooth boundary, over which the energy we are interested in, i.e. their perimeter, is classically defined. However a sequence of smooth sets with uniform controlled perimeter might converge to nonsmooth sets (for instance one can easily form corner, cusps, and other type of singularities) and we therefore would like to enlarge this class. We then take a much larger class, that of all measurable sets, with a topology in which the good objects are dense and we extend the energy to be the lower semicontinuous envelope.  Interestingly Caccioppoli's approach was initially dismissed by his contemporaries (cf. the reviews by L.C. Young of the aforementioned papers) because he was not able to relate his abstract definition to any concrete notion of perimeter in a measure-theoretic sense. In the early fifties De Giorgi took up Caccioppoli's approach and proved, in his celebrated works on the isoperimetric property of the sphere (cf. \cite{DG1,DG2,DG3}), that:
\begin{itemize}
\item the class of sets with finite perimeter in the sense of Caccioppoli is compact, under a uniform bound on their perimeter;
\item the perimeter has a precise measure-theoretic interpretation, i.e. if the set $\Omega$ has a finite perimeter one can introduce a suitable notion of (oriented) measure-theoretic boundary which turns out to be rectifiable and whose Hausdorff measure is indeed the perimeter of $\Omega$. 
\end{itemize}
De Giorgi also reformulated the theory of sets of finite perimeters through a useful duality: if correctly interpreted, the usual divergence theorem holds for them, and the boundary integral in the formulation is in fact a classical integral, in the sense of measure theory, over the measure-theoretic boundary. For open sets with smooth boundaries the ``measure-theoretic'' one coincides with the topological one. An interesting byproduct (not at all obvious from the definition) is that the perimeter as defined by Caccioppoli is in fact the classical surface area of the topological boundary when the latter is smooth.

Thus, the oriented ``generalized'' boundaries of Caccioppoli and De Giorgi act as linear functionals on vector fields. In the celebrated theory developed later by Federer and Fleming (cf. \cite{FF}) these are particular instances of ``integral currents'', which act on general forms (and hence can have arbitrary codimension). Like De Giorgi's theory of Caccioppoli sets, the theory of integral currents of Federer and Fleming can also be seen as a suitable variational completion: after introducing an appropriate class of good objects (in this case integral smooth chains, which are formal linear combinations, with integer coefficients, of smooth oriented submanifolds with smooth boundaries, the more general objects, namely the integral currents, can be characterized as the limits, in an appropriate  weak topology, of sequences of those good objects, under uniform bound on their volume and on the volume of their boundaries. Like in the case of De Giorgi's theory of Caccioppoli sets, integral currents can be represented, in a suitable measure theoretic sense, as integration over ``oriented'' rectifiable sets. 

While the duality with differential forms limits the choice of coefficient groups in the formal linear combinations to integer and real coefficients (or anyway to subgroups of the reals), the ``completion point of view'' allows to choose other ``coefficient groups'' (endowed with an appropriate norm, so that we can make sense of the notion of ``mass''), cf. the foundational paper of \cite{Fleming} for the case of finite groups. Notable choices are the so called ``flat chains mod $p$'' (which with a slight abuse of terminology we will call currents mod $p$). In the latter case $p$ is a positive integer larger than $1$ and the coefficient group is $\mathbb Z_p = \mathbb Z/ (p\mathbb Z)$ (for an element $[q]\in \mathbb Z_p$, endowed with the usual norm
\[
|[q]|\min \{|q-kp|: k\in \mathbb Z\}\, .
\] 
In this note the coefficients group will always be either $\mathbb Z$, or $\mathbb Z_p$. Note in particular that in both cases the norm will always take integer values, a fact which will play a fundamental role in our discussions.

In all these instances we have a framework where we can apply the direct methods of the calculus of variations. In particular:
\begin{itemize}
\item the concept of boundary comes naturally either from the duality with differential forms, or from the closure procedure;
\item  the underlying space of  generalized objects is closed;
\item the generalized area functional (often called the mass) is lower semicontinuous and its sublevel sets are compact (if we assume that the boundary of our generalized surfaces is a fixed given one).
\end{itemize}
In particular the Plateau problem in the above frameworks has a very elegant existence theory. 

\subsection{Varifolds and the calculus variations ``in the large''} One notable drawback of the functional analytic frameworks outlined above is that the mass is not continuous for the natural convergence in the underlying spaces. Continuity along a sequence might be lost because of two mechanisms:
\begin{itemize}
\item High frequency oscillations: for instance the graphs of the functions $\frac{1}{k} \sin kx$ in the two-dimensional plane have locally bounded length and they converge, in the sense of integral currents, to the straight line. It is however easy to see that the total length of any segment in the limiting line is strictly less than the limit of the corresponding approximations.
\item Cancellation: a line in $\mathbb R^2$ can be given two distinct orientations, thereby defining two different integral currents. However their sum is $0$. If we approximate the two different oriented lines with a sequence of two shifted oriented lines with disjoing supports, then we get a sequence of integral currents with masses uniformly bounded from below which converges to the trivial current.     
\end{itemize} 
We are concerned mostly with the second, since the approximating sequence is a sequence of minimal surfaces. Indeed we can reasonably expect that a sequence of critical surfaces will not exhibit the oscillatory behavior of the first example (this fact was in fact proved by Allard as a byproduct his famous regularity theory, see below). On the other hand the criticality assumption does not rule out the second example, which is therefore ``particularly bad'' because it shows that in the space of currents we cannot expect any reasonable type of ``Palais-Smale'' property.

A way to remedy this loss of continuity is to introduce the notion of varifold, which is just a positive measure on the Grassman space of $m$-dimensional unoriented $m$-planes in the tangent bundle of the Euclidean space (or more generally of a Riemannian manifold). For a smooth (not necessarily oriented) surface $\Sigma$, the corresponding varifold is given by 
\begin{equation}\label{e:varifolds}
\delta_{T_x \Sigma} \otimes d\textrm{Vol}_\Sigma
\end{equation} 
General varifolds were introduced by L. C. Young (cf. \cite{Young}), while in the context of minimal surfaces Almgren introduced them precisely in order to tackle general existence problems for {\em critical} points of the area functional, cf. \cite{Alm65}. A particularly useful subclass of varifolds is that of integral varifolds, which satisfy a structure as in \eqref{e:varifolds} where $d\textrm{Vol}_\Sigma$ is substituted by the Hausdorff $k$-dimensional measure restricted on a general $m$-dimensional rectifiable set $R$ (with an integer valued weight), and $\delta_{T_x\Sigma}$ is substituted by $\delta_{T_x R}$, where $T_x R$ is a natural measure-theoretic generalization, to rectifiable sets, of the tangent to a smooth surface.  

In his notable monograph \cite{Pitts}, based on some groundbreaking ideas of Almgren \cite{Alm65}, Pitts developed a quite powerful variational theory for finding generalized critical points of the area functional. In codimension 1, i.e. in the case of hypersurfaces, the theory of Almgren and Pitts has found striking geometric applications in the works of A. Neves and F. Coda Marques, cf. \cite{MN1,MN2,MN3}. These results have spurred a number of interesting works in the area and Pitts' existence theory has been revisited in several different ways, see for instance \cite{CM,CD,DT,Ket,ZZ}. 

Varifolds can be naturally deformed using one parameter family of diffeomorphisms and this allows to introduce a rather natural notion of $k$-th variation of the varifold along smooth vector fields. Of particular relevance are then
\begin{itemize}
\item {\em stationary} varifold, i.e. varifolds for which the first variation vanishes along any vector field, 
\item and {\em stable} varifolds, i.e. stationary varifolds for which the second variation is nonnegative along any vector field. 
\end{itemize}
Since all the objects encountered above in the existence theories for the Plateau problem naturally induce corresponding varifolds, all the minimizers in the various senses given above are in fact stable varifolds. 

\section{Monotonicity formula and tangent cones}

One simple and very powerful tool in the regularity theory for minimal submanifolds is the monotonicity formula. In order to gain an intuition about it, consider a smooth $m$-dimensional surface $\Sigma\subset \mathbb R^{m+n}$ which minimizes the volume in some suitable class of comparison surfaces and fix an ``interior'' point $p\in \Sigma$. Considerthen  a ball $\bB_r (p)$ which does not  intersect $\partial \Sigma$. We wish to compare the volume of $\Sigma \cap \bB_r (p)$ to the volume of the cone $\Sigma_c$ with vertex $p$ and base $\Sigma \cap \partial \bB_r (p)$. By Sard's lemma we can assume that $\Sigma$ intersects $\partial \bB_r (p)$ transversally. Note that our comparison surface is somewhat singular, because of the vertex singularity of the cone and the discontinuity in the tangents that might be introduced by cutting $\Sigma \cap \bB_r (p)$ out and replacing it with $\Gamma$. On the other hand it is also simple to see that $\Sigma_c \cup (\Sigma\setminus \bB_r (p))$ can be obtained as limit of deformations of $\Sigma$ by smooth isotopies of the ambient space: in particular it is a good comparison surface in pretty much all the variational frameworks considered so far.

The minimizing property of $\Sigma$ implies then that 
\[
\textrm{Vol}^m (\Sigma \cap \bB_r (p)) \leq \textrm{Vol}^m (\Sigma_c) = \frac{r}{m} \textrm{Vol}^{m-1} (\Sigma \cap \partial \bB_r (p))\, ,
\]
which in turn (considering that $\textrm{Vol}^{m-1} (\Sigma \cap \partial \bB_r (p))\leq \frac{d}{dr} \textrm{Vol}^m (\Sigma \cap \bB_r (p))$) gives 
\begin{equation}\label{e:monotonicity}
\frac{d}{dr} \frac{\textrm{Vol}^m (\Sigma \cap \bB_r (p))}{r^m}\geq 0 \, .
\end{equation}
The latter is the classical monotonicity formula for minimal submanifolds. It is very robust, in the sense that
\begin{itemize}
\item[(a)] It can be derived for critical points by using stationarity with respect to some specific radial deformations. In particular it holds for stationary varifolds (see \cite{Allard72}).
\item[(b)] Allowing for suitable multiplicating factors like $e^{Cr}$, the formula holds for much more general objects, in particular for stationary varifolds in smooth Riemannian manifolds (cf. again \cite{Allard72}).
\item[(c)] A suitable version of the formula can be derived at boundary points too, under the assumption that the boundary $\partial \Sigma$ is smooth enough (cf. \cite{Allard75}). An intuition for this can be gained through the following observation: if the boundary $\partial\Sigma$ were an affine subspace passing through $p$, then the competitor surface $\Sigma_c \cup (\Sigma\setminus \bB_r (p))$ has the same boundary, namely $\partial \Sigma$. 
\item[(d)] Perhaps most importantly, a more refined version of the arguments leading to \eqref{e:monotonicity} shows that the equality case in holds if and only if $\Sigma_c$ coincides with $\Gamma$, i.e. if $\Sigma$ itself is a cone with vertex $p$.
\end{itemize}
For further reference we will call {\em density of $\Sigma$ at $p$} (denoted by $\Theta (\Sigma, p)$  the limit of the ``mass ratio''
\[
\lim_{r\downarrow 0}  \frac{\textrm{Vol}^m (\Sigma \cap \bB_r (p))}{\omega_m r^m}\, ,
\] 
where $\omega_m$ is the volume of the $m$-dimensional disk. Obviously the density is not particularly interesting for smooth $\Sigma$'s as it will be $1$ at every interior point and $\frac{1}{2}$ at every boundary point (or $k$ at every interior point and $\frac{k}{2}$ at every boundary point, if we entertain the possibility of allowing for multiplicities in the normed groups $\mathbb Z$ and $\mathbb Z_p$, or if we consider integral varifolds). However, due to (a), the density exists for any ``generalized minimal surface'' encountered in the previous section and this is a very nontrivial information, given that the latter might be singular. Another interesting byproduct of the monotonicity formula is that the density is nowhere smaller than $1$ at interior points, which in turn implies that some suitable definition of ``support'' of the generalized minimal surface is a closed rectifiable set of locally finite Hausdorff measure.

\subsection{Tangent cones}
Fact (d) above is maybe the most relevant, as it is the starting point for a fruitful fundamental concept in minimal surface theory. Let us fix a point $p$ in (the support of) our generalized minimal surface $\Sigma$ (and, in light of (c) above, we might even fix it at the boundary as long as the latter is sufficiently smooth). For every radius $r$ consider then the translated and rescaled surface
\[
\Sigma_{p,r} := \frac{\Sigma-p}{r} = \{y: p+ry \in \Sigma\}\, .
\] 
The volume of the surface in $\bB_R (0)$ is then uniformly bounded for every fixed $R$, independently of the parameter $r$. Again, while this is not particularly exciting for a smooth $\Sigma$, it is a highly nontrivial amount information for the generalized minimal surfaces, which are potentially singular at $p$. Given the uniform bound and the compactness properties available for all the generalized minimal surfaces introduced thus far, up to subsequences we can assume that $\Sigma_{p, r}$ converges to a generalized minimal surface in the same class, which for convenience we will denote by $\Sigma_c$. 

Assuming convergence of the volume (which is in fact correct for objects like varifolds because of their definition, while it is a property shared by minimizers out of variational arguments, in any of the classes described thus far) the mass ratio $R^{-m} \textrm{Vol}^m (\Sigma_c \cap \bB_R)$ is constant in $R$, and hence by point (d) it is a cone. In the literature $\Sigma_c$ is called {\em a tangent cone to $\Sigma$ at $p$}. Note that we are speaking about {\em a} tangent cone: the uniqueness of this object, i.e. the independence of it from the subsequence $r_k\downarrow 0$, is a widely open problem, even though several fundamental result have been proved in the past (see Section \ref{s:tangent_cones} below). 

At a regular point $p$, i.e. a point in a neighborhood of which the generalized minimal surface $\Sigma$ is smooth, the tangent cone $\Sigma_c$ is of course unique and it is given by the tangent space to $\Sigma$ at $p$ (counted with the appropriate multiplicity, depending upon the chosen variational framework), or half of the tangent space if $p$ is a boundary point. A later section will examine under which assumption the latter conclusion is correct. 

At any rate, even in the possible presence of singularities, we have gained a great deal of new information about $\Sigma_c$ compared to $\Sigma$: $\Sigma_c$ is a ``global minimal surface'' (it has no boundary if $p$ is in the interior, or its boundary is affine if $p\in \partial \Sigma$) and moreover it is conical. In particular its spherical cross section carries all the information about $\Sigma_c$, even  when $\Sigma_c$ is singular: at all effects $\Sigma_c$ must be less complex than $\Sigma$, i.e. $\Sigma_c$ has ``lost one dimension''.  

\section{Invariant spaces and strata}\label{s:strata}

 For simplicity in what follows we will focus on tangent cones $\Sigma_c$ at interior points $p$, even though a variant of the following discussion applies to boundary tangent cones as well. Since $p$ is not a boundary point, $\Sigma_c$ has no boundary. Moreover two simple corollaries of the monotonicity formula are that:
\begin{itemize}
\item $\Theta (\Sigma_c, 0) \geq \Theta (\Sigma_c, q)$ for every $q$,
\item and if the quality holds at some $q\neq 0$,  then $\Sigma_c$ ``splits off a line'', i.e. it is invariant under translations in the direction $q$.
\end{itemize}
The latter property is the starting point of Federer in his celebrated ``dimension reduction argument'' (cf. \cite{Federer,Federer70}), which we will illustrate below. Here we want to present Almgren's stratification theory, which is a far-reaching generalization of Federer's original idea. 

First of all, it follows from the above consideration that the set 
\[
V= \{q\in \Sigma_c: \Theta (\Sigma_c, q) = \Theta (\Sigma_c, 0)\}
\] 
is a linear subspace of $\mathbb R^{m+n}$. If $V$ has the same dimension as $\Sigma_c$, then in fact $\Sigma_c$ coincides with $V$ (counted with the correct multiplicity and, in some cases, given the correct orientation). Otherwise, assuming that $k = \textrm{dim} (V)$, $\Sigma_c$ is the product of $V$ and a minimal cone $\Sigma_0$ in the orthogonal complement of $V^\perp$, which is not invariant by any translation. This is a great deal of information and in several cases implies severe restrictions upon $k$. For instance, for area minimizing integral currents it can be easily checked that $k\leq m-2$ in general, while in the particular case of codimension $1$ the celebrated paper of Simons on stable minimal hypercones (cf. \cite{Federer}) implies that $k\leq m-7$ (again this will be discussed further below)! In \cite{Alm00} Almgren coined the term {\em building dimension of the cone $\Sigma_c$} to identify the nonnegative integer $k$ and introduced a stratification of the interior points $p\in \Sigma$ according to the maximal building dimension of its tangent cones. In particular the stratum $\mathscr{S}_k$ is the set of (interior) points $p\in \Sigma$ such that the building dimension of any tangent cone to $\Sigma$ at $p$ is at most $k$. Almgren's fundamental discovery is the following

\begin{theorem}\label{t:Almg-strat}
For a stationary integral varifold $\Sigma$ the stratum $\mathscr{S}_k$ is a closed set of Hausdorff dimension at most $k$.
\end{theorem}

Almgren's approach is very general and can be applied to a variety of different context. For a framework which is very flexible and covers a white range of applications, see \cite{White97}.
Almost four decades after the work of Almgren, groundbreaking ideas allowed Naber and Valtorta to improve massively upon Almgren's original theorem, showing that 
(cf. \cite{NV1,NV2})

\begin{theorem}\label{t:NV}
For a stationary integral varifold $\Sigma$ the stratum $\mathscr{S}_k$ is $k$-rectifiable, i.e. it can be covered, up to a set of $\mathcal{H}^k$-measure zero, with countably many $C^1$ submanifolds of dimension $k$.
\end{theorem}

Theorem \ref{t:NV}, which was predated by pioneering works of Simon (cf. \cite{Simon93,Simon96}) covering some particular cases (most notably the stratum $\mathscr{S}_{m-7}$ for area-minimizing integral currents, see below for more details), builds upon a new sophisticated version of Reifenberg's topological disk theorem combined with a clever use of the remainder in the monotonicity formula. The ideas are quite general and can be applied to other contexts. 

\section{Interior $\varepsilon$-regularity at multiplicity $1$ points}

Following the above terminology, two things are obvious: the stratum $\mathscr{S}_m$ coincides with the whole support of the $m$-dimensional generalized minimal surface and the stratum $\mathscr{S}_{m-1}$ consists necessarily of singular points. A point $p\in \mathscr{S}_m\setminus \mathscr{S}_{m-1}$ is clearly a good candidate for being a regular point, since we know that at least one tangent cone to $\Sigma$ at $p$ is in fact a plane (counted with its multiplicity). However, a famous theorem by Federer shows that the existence of a ``flat tangent'' does not guarantee the regularity of the point. Indeed, based on a classical theorem of Wirtinger in K\"ahler geometry, Federer proved (cf. \cite{Federer}) that

\begin{theorem}\label{t:Federer}
Any holomoprhic subvariety $\Sigma$ of complex dimension $k$ in $\mathbb C^n$ induces an area-minimizing integral current of dimension $2k$ in $\mathbb R^{2n}$.
\end{theorem}

It can be readily checked that the holomorphic curve 
\begin{equation}\label{e:example}
\Sigma = \{z^2=w^3: (z,w)\in \mathbb C^2\}\, .
\end{equation} 
gives then an example of an area-minimizing integral current of dimension $2$ in $\mathbb R^4$ for which $0\in \mathscr{S}^2\setminus \mathscr{S}^1$ is a singular point. One crucial fact is however that the flat tangent at $0$ is a $2$-dimensional plane (i.e. the complex line $\{z=0\}$) but counted with {\em multiplicity $2$}. A celebrated theorem of Allard (cf. \cite{Allard72}), extensively used in the literature, shows that the naive expectation ``flat tangent cone $\iff$ regular point'' is indeed correct if the flat tangent cone has multiplicity $1$.

\begin{theorem}\label{t:Allard}
If a stationary integral varifold $\Sigma$ is sufficiently close in $\bB_{2r} (p)$ to a plane (counted with multiplicity $1$) in the weak topology, then in $\bB_r (p)$ it is a smooth graph over that plane. Moreover, at any interior point $p$ where the density of $\Sigma$ is $1$, such a plane always exists for a sufficiently small $r$.
\end{theorem}

Among the various objects examined in this note, there are three situations where it is relatively simple to see apriori that our generalized minimal surface $\Sigma$ will not ``pick higher multiplicity'' at flat points:
\begin{itemize}
\item[(a)] $\Sigma$ is a portion of the boundary of some Caccioppoli set;
\item[(b)] $\Sigma$ is a solution of the Plateau problem in one of the set-theoretic senses described in Section \ref{ss:set-theoretic};
\item[(c)] $\Sigma$ is an area-minimizing current mod $2$ or an area-minimizing current mod $3$. 
\end{itemize}
In fact Theorem \ref{t:Allard} was realized independently by De Giorgi and Reifenberg, in \cite{DG4} and \cite{Reif2,Reif3}, respectively in the particular cases of (a) and (b) (this is literally correct for De Giorgi, while in reality Reifenberg in \cite{Reif2,Reif3} dealt with the only set-theoretic solutions of the Plateau problem known at his time, which were the ones he himself introduced in \cite{Reif1}; it must also be noticed that De Giorgi's monograph appeared 3 years before Reifenberg's paper, but it was probably not yet widely known when Reifenberg wrote his papers \cite{Reif2,Reif3}). The two pioneering approaches are rather different, but they both rely on the fact that the ``linearization'' of a minimal surface, understood as a graph over his tangent plane, is harmonic (in fact it would be more correct to say that the linearization of the minimal surface equation is the Laplace equation, or that, at the level of the energies, the Dirichlet energy is the second order Taylor expansion of the area functional). 

Reifenberg used harmonic competitors to estimate how much an area-minimizing surface deviates from being conical if it is close to a plane, and derived his famous ``epiperimetric inequality'', which can be thought as a quantitative improvement of the cone-comparison outlined above to prove the monotonicity formula. De Giorgi used a linearization technique which has a more PDE flavor, and which was generalized afterwards by Almgren in any codimension and for much more general energy functionals, cf. \cite{Alm65}. Both approaches exploit in a substantial way the minimizing property of the surfaces in question. Allard's proof of Theorem \ref{t:Allard}, while still based on the intuition that harmonic functions provide a good approximation for minimal graphs, deviates drastically from both of them, having to deal with stationary objects. But ultimately it is fair to say that Allard's approach borrows much more substantially from the works of De Giorgi and Almgren, than from that of Reifenberg. 

It is worth spending some words on why all the approaches mentioned above for the $\varepsilon$-regularity theory fail at the origin in the example \eqref{e:example}: no matter how small is the scale that we look at, it is not possible to approximate efficiently \eqref{e:example} around the origin with the graph of a single-valued function. Of course, before knowing Theorem \ref{t:Allard} we also do not know that, under the corresponding assumptions, a generalized surface is graphical over the approximating plane: however a crucial point in Allard's proof of Theorem \ref{t:Allard} is that, before proving any regularity, he was able to produce a graphical approximation which covers {\em most} of the support of the generalized minimal surface. In contrast, no matter how small the $r$ is, a single valued graph will cover no more than half of $\Sigma \cap \bB_r (p)$ when $\Sigma$ is given by \eqref{e:example}
 
The assumption on the multiplicity of the varifold severely limits the effectiveness of Theorem \ref{t:Allard} in bounding the size of the singular set for stationary integral varifolds. In fact, it would be natural to expect that singular points with a flat tangent cone form anyway a set of relatively modest size: according to the known examples, its dimension is likely $m-2$. The latter is less than the dimension of $\mathscr{S}_{m-1}$ and so one could reasonably conjecture that the singular set of a stationary integral varifold has dimension at most $m-1$. On the other hand so far the best that we can conclude is still a corollary of Theorem \ref{t:Allard} noted by Allard in \cite{Allard72} almost 50 years ago.

\begin{corollary}\label{c:Allard}
Let $\Sigma$ be a stationary $m$-dimensional varifold in $U\subset \mathbb R^{m+n}$. Then the singular set of $\Sigma$ is a closed subset which has empty (relative) interior. 
\end{corollary}

\section{Boundary $\varepsilon$-regularity at multiplicity $\frac{1}{2}$ points}

In his second groundbreaking work \cite{Allard75} Allard proved a statement parallel to Theorem \ref{t:Allard} at boundary points. The following is an informal description of his main ``boundary regularity'' theorem.

\begin{theorem}\label{t:Allb}
Assume $\Sigma$ is an $m$-dimensional integral varifold in some open set $U\subset \mathbb R^{m+n}$, which is stationary for variations which keep fixed a smooth $m-1$-dimensional submanifold $\Gamma$. Then the following conclusions hold:
\begin{itemize}
\item[(a)] If $p\in \Gamma$ belongs to the support of the varifold, then $\Theta (p, \Sigma) \geq \frac{1}{2}$;
\item[(b)] If in $\bB_{2r} (p)$ the varifold is sufficiently close, in the weak topology, to a single copy of half of an $m$-dimensional plane $\pi$, then in $\bB_r (p)$ it is a $C^1$ graph over a suitable portion of $\pi$;
\item[(c)] If $\Theta (p, \Sigma)=\frac{1}{2}$, then the assumption of (b) (and hence the corresponding conclusion) holds for a sufficiently small $r$.
\end{itemize}
\end{theorem}

For this boundary version as well, the overall intuition is that $V$ is, in first approximation, very well approximated by the graph of a (single-valued) harmonic function. 

While in the rest of this note I will touch upon interior regularity results for many different notions of generalized minimal surfaces, concerning boundary regularity I will only focus on the case of area minimizing integral currents. This is also due to the fact that there are not many other cases studied in the literature. Aside from Allard's general theorem (i.e. Theorem \ref{t:Allb} stated above), the author is only aware of:
\begin{itemize}
\item the work \cite{LM} (cf. also \cite{Morgan}), which contains a conjectural list of boundary tangent cones for set-theoretic $2$-dimensional solutions of the Plateau problem;
\item the recent work of David \cite{David4}, which, for $2$-dimensional sliding minimizers, generalizes the conclusion of Theorem \ref{t:Allb} to the union of two half planes, and possible additional transverse cones as in the classical theorem of Taylor in the interior (cf. Theorem \ref{t:Taylor});
\item an argument by White which shows how to gain curvature estimates for stable minimal hypercurrents at the boundary, under some convexity assumption (cf. \cite[Section 6.4]{DR}). 
\end{itemize} 

\section{Interior regularity theory: minimizing integral hypercurrents} 

Even though Allard's Theorem \ref{t:Allard} needs the multiplicity $1$ assumption, the latter might be dropped in the case of integral area-minimizing currents of codimension $1$ (which for simplicity we will call hypercurrents from now on). The key point is that area-minimizing integral hypercurrents $\Sigma$ can be locally decomposed into the sum of area-minimizing boundaries of Caccioppoli sets (this is a consequence of  the Coarea Formula, see for instance \cite{Simon-book}). If in $\bB_{2r} (p)$ the original current $\Sigma$ is close to a multiple $Q$ of a hyperplane $\pi$, each of these boundaries is then close to a multiplicity $1$ copy of $\pi$. We can then apply Allard's theorem to prove that each of them is a $C^1$ graph in $\bB_r (p)$, obtaining what can be called (cf. for instance \cite{SS}) a ``sheeting theorem'' for $\Sigma\cap \bB_r (p)$. However each of these sheets must be ordered (for minimizing reasons they cannot cross) and they must touch at the point $p$: the maximum principle (each of these graphs is a solution of the minimal surface equation) then implies that they collapse all into a single smooth surface counted with the appropriate multiplicity, which must be $Q$. 

This argument rules out that an example like \eqref{e:example} could exist for integral area-minimizing {\em hyper}currents. We are therefore in the luckiest of situations where we can infer that a single flat tangent cone at $p$ is indeed a necessary and sufficient condition for regularity at $p$. If we introduce the notation $\Sing_i (\Sigma)$ for the interior singularities of $\Sigma$, when $\Sigma$ is an $m$-dimensional area-minimizing integral current in $\mathbb R^{m+1}$ (or more generally in a complete smooth Riemannian manifold of dimension $m+1$) we infer $\Sing_i (\Sigma)\subset \mathscr{S}_{m-1}$. Consider, however, that the existence of a point $p\in \mathscr{S}_{m-1}\setminus \mathscr{S}_{m-2}$ implies the existence of a singular $1$-dimensional area-minimizing cone in $\mathbb R^2$, and it is rather elementary to see that the latter cones do not exist, namely that $\Sing_i (\Sigma) \subset \mathscr{S}^{m-2}$. Of course we can now wonder whether for some $\Sigma$ the set $\mathscr{S}_{m-2}\setminus \mathscr{S}_{m-3}$ is nonempty, which is equivalent to the existence of an area-minimizing $2$-dimensional cone in $\mathbb R^3$ which is not a plane (i.e. it is singular at the origin). In \cite{Federer70} Federer introduced his well known reduction argument, which could be formalized as follows.

\begin{theorem}\label{t:Federer2}
Let $m$ be the smallest integer with the property that there is an $m$-dimensional area-minimizing integral current $\Sigma_0$ in $\mathbb R^{m+1}$ which is a nonplanar cone with vertex at the origin. Then  $\Sigma_0$ is everywhere regular except at the origin.
\end{theorem}

It was also realized by De Giorgi in \cite{DG5} that the well known Bernstein problem, i.e. whether a complete minimal graph over $\mathbb R^{m+1}$ must be affine, would also be implied by the nonexistence of nonplanar area-minimizing oriented hypercones in $\mathbb R^{m+1}$. After progress by Fleming, De Giorgi, and Almgren (cf. \cite{Fleming62,DG5,Alm66}), Simons in \cite{Simons} proved his famous result about stable minimal hypercones, namely

\begin{theorem}\label{t:Simons}
If $m\leq 6$ and $\Gamma_0 \subset \partial \bB_1\subset \mathbb R^{m+1}$ is a smooth connected submanifold of dimension $m-1$, such that the cone $\Sigma_0$ with base $\Gamma_0$ and vertex $0$ is a stable varifold, then $\Gamma_0$ is a great sphere (i.e. $\Sigma_0$ is planar). On the other hand 
\begin{equation}\label{e:Simons}
\Sigma_s := \{x_1^2+x_2^2+x_3^2+x_4^2 = x_5^2+x_6^2+x_7^2+x_8^2\} \subset \mathbb R^8
\end{equation}
is a nonplanar, oriented, stable singular cone of dimension $7$.
\end{theorem}

Since area-minimizing currents are automatically stable varifolds, in combination with Federer's reduction argument the first part of Theorem \ref{t:Simons} implies that $\Sing_i (\Sigma)\subset \mathscr{S}_{m-7}$ for any $m$-dimensional area-minimizing integral hypercurrents. In particular for $m\leq 6$, $\Sigma$ is a regular hypersurface (in the interior), while, for $m\geq 7$, $\Sing_i (\Sigma)$ has dimension at most $m-7$. In fact, by Theorem \ref{t:NV} we can conclude that $\Sing_i (\Sigma)$ is $m-7$-rectifiable. The latter conclusion was first reached by Simon in his pioneering work \cite{Simon93}. However, compared to Simon's techniques, the approach by Naber and Valtorta (cf. \cite{NV2}) allows to prove the stronger conclusion that

\begin{theorem}\label{t:NV-2}
Let $\Sigma$ be an $m$-dimensional area-minimizing integral current in $\mathbb R^{m+1}$. Then $\Sing
_i (\Sigma)$ has locally finite Hausdorff $m-7$-dimensional measure, and it is $m-7$-rectifiable. 
\end{theorem} 

In their famous work \cite{BDG} Bombieri, De Giorgi and Giusti completed the solution of the Bernstein problem showing that indeed the Simons cone \eqref{e:Simons} is an area minimizing integral current of dimension $7$, and that in addition there is a nonaffine global solution $u: \mathbb R^8\to \mathbb R$ of the minimal surface equation.

\section{Interior regularity theory: minimal sets}

As already mentioned, the phenomenon of ``picking higher multiplicity at flat points'' is absent in the solutions of the Plateau problem that fall in the ``set-theoretic'' approach. This was pioneered by Reifenberg in \cite{Reif2,Reif3}, who proved that his $m$-dimensional solutions of the Plateau problem are always real analytic except for a closed $\mathcal{H}^m$-null set. A much more general statement, valid in a variety of contexts and also for a vast class of elliptic energies was proved by Almgren in \cite{Alm76}. 

Following the Remarks of Section \ref{s:strata} we conclude that the singular set of an $m$-dimensional set-theoretic solution of the Plateau's problem is necessarily contained in $\mathscr{S}^{m-1}$. While Theorem \ref{t:NV} implies that $\mathscr{S}^{m-1}$ is rectifiable, much more can actually be said in the codimension $1$ case. First of all, for $2$-dimensional minimizing sets in $\mathbb R^3$ Taylor in \cite{Taylor76} proved the following complete structure theorem.

\begin{theorem}\label{t:Taylor}
Let $\Sigma$ be a $2$-dimensional set which minimizes the area in the sense of Almgren. Then:
\begin{itemize}
\item[(a)] $\mathscr{S}_1 (\Sigma)\setminus \mathscr{S}_2 (\Sigma)$ is the (locally finite) union of $C^{1,\alpha}$ arcs and for each $p\in \mathscr{S}_1 (\Sigma)$ there is a neighborhood $U$ of $p$ in which $\Sigma$ is the union of three classical minimal surfaces meeting in $\mathscr{S}_1 (\Sigma) \cap U$ at $120$ degrees;
\item[(b)] $\mathscr{S}_0 (\Sigma)$ consists of isolated points and for each $p\in \mathscr{S}_0 (\Sigma)$ there is a neighborhood $U$ in which $\Sigma$ is diffeomorphic to the cone over a regular tetrahedron.  
\end{itemize}
\end{theorem}

The same conclusion as in part (a) of the above remarkable theorem is in fact valid for the stratum $\mathscr{S}_{m-1}\setminus \mathscr{S}_{m-2}$ of  $m$-dimensional area-minimizing sets in $\mathbb R^{m+1}$ and can be inferred from Simon's theory on the uniqueness of multiplicity $1$ cylindrical cones, cf. \cite{Simon93b}. Part (b) can also be genealized to a similar statement for $m$-dimensional area-minimizing sets of $\mathbb R^{m+1}$, implying in particular that $\mathscr{S}_2\setminus \mathscr{S}_{m-3}$ is an $m-2$-dimensional submanifold. This generalization was announced by White in \cite{White86} and a proof has been recently published by Colombo, Edelen, and Spolaor in \cite{COS}, as a corollary of a more general result. The main Theorem in \cite{COS} also implies that the stratum $\mathscr{S}_{m-3}$ has finite $\mathcal{H}^{m-3}$ measure. 

\section{Interior regularity theory: stable hypersurfaces and stable hypervarifolds}

In \cite{SSY} Schoen, Simon, and Yau realized that Simons' Theorem on stable minimal hypercones could be recast in a suitable apriori estimate for the curvature of stable minimal surfaces. More precisely, combining Simons' inequality with techniques from elliptic PDEs they were able to prove the following groundbreaking theorem.

\begin{theorem}\label{t:SSY}
Let $\Sigma$ be a smooth minimal hypersurface in $U\subset \mathbb R^{m+1}$ with $m\leq 5$. Then for every $V\subset\subset U$ there is a constant $C$ which depends on $U,V$, and $\mathcal{H}^m (\Sigma)$ such that the Hilbert-Schmidt norm of the second fundamental form $A$ of $\Sigma$ is bounded by $C$ at every point of $\Sigma\cap V$.
\end{theorem}

In their subsequent work \cite{SS} Schoen and Simon were able to cover the case $m=6$ of the above statement and also to give a ``GMT regularity theory'' counterpart of the Schoen-Simon-Yau estimates. More precisely they were able to prove

\begin{theorem}\label{t:SS}
Assume $\Sigma$ is a stable $m$-dimensional varifold in $U\subset \mathbb R^{m+1}$ with the property that\footnote{In their paper Schoen and Simon assume the stronger property that $\mathcal{H}^{m-2} (\Sing_i (\Sigma))=0$, but it is well known by the experts that their arguments apply if the Hausdorff measure is in fact finite.} $\mathcal{H}^{m-2} (\Sing_i (\Sigma)) < \infty$. Then $\Sing_i (\Sigma) \subset \mathscr{S}_{m-7}$.
\end{theorem}

It has been recently shown by Simon, see \cite{Simon21a,Simon21b}, that the subsequent conclusion that $\Sing_i (\Sigma)$ is $m-7$-rectifiable is optimal, in the sense that there are stable $m$-dimensional varifolds in $m+1$ smooth Riemannian manifolds whose singular sets are closed sets of arbitrary Hausdorff dimension $\alpha \leq m-7$. On the other hand the assumption $\mathcal{H}^{m-2} (\Sing_i (\Sigma))=0$ is not at all optimal. Based on the examples known thus far, one could expect that for a general stable hypervarifold the top stratum $\mathscr{S}_{m-1}\setminus \mathscr{S}_{m-2}$ is a $C^{1,\alpha}$ $m-1$-dimensional submanifold and that if the latter is empty, then $\Sing_i (\Sigma)$ is contained in $\mathscr{S}_{m-7}$. A notable theorem in this direction, in particular covering the second conclusion, has been achieved by Wickramasekera in his deep regularity theory of stable hypervarifolds. The main conclusion of his paper \cite{Wic} is the following

\begin{theorem}\label{t:Wick}
Assume $\Sigma$ is a stable $m$-dimensional varifold in a connected set $U\subset \mathbb R^{m+1}$. Then:
\begin{itemize}
\item either $\Sing_i (\Sigma)$ contains a point $p$ in a neighborhood of which $\Sigma$ consists of a finite number of smooth minimal hypersurfaces meeting at a commong $C^{1,\alpha}$ $m-1$-dimensional boundary (which in particular is a nonempty subset of $\mathscr{S}_{m-1}\setminus \mathscr{S}_{m-2}$),
\item or otherwise $\Sing_i (\Sigma) \subset \mathscr{S}_{m-7}$.
\end{itemize}
\end{theorem}

While the latter is a remarkable achievement, for general stable hypervarifolds the best nonconditional regularity result is still the one that can be concluded from the sole condition of stationarity through Allard's work, namely Corollary \ref{c:Allard}.

\section{Interior regularity theory: minimizing integral currents in higher codimension}

As already witnessed in Example \eqref{e:example} the regularity theory for area-minimizing integral currents in codimension larger than $1$ differs dramatically from the regularity theory for hypercurrents, since there are singular points which belong to $\mathscr{S}_m\setminus \mathscr{S}_{m-1}$, and which from now on we will call ``flat singular points''. The major problem of how to give a suitable dimension bound for ``flat singular points'' was finally conquered by Almgren in a titanic effort, which resulted in a famous 1728 pages preprint in the early eighties (cf. \cite{Alm83}), published pusthumously thanks to the editorial work of Scheffer and Taylor in \cite{Alm00}. Almgren's monograph achieves the optimal dimension bound for area-minimizing integral currents in any dimension and codimension.

\begin{theorem}\label{t:Almgren}
Let $\Sigma$ be an area-minizing integral current of dimension $m$ in $\mathbb R^{m+n}$. Then the Hausdorff dimension of the set of interior flat singular points is at most $m-2$, while the stratum $\mathscr{S}_{m-2}\setminus \mathscr{S}_{m-1}$ is empty. In particular $\textrm{dim}_H (\Sing_i (\Sigma)) \leq m-2$.
\end{theorem}

Almgren invented several tools to prove Theorem \ref{t:Almgren}. In particular:
\begin{itemize}
\item[(i)] he introduced an entire new concept of ``multivalued functions minimizing the Dirichlet energy'' in order to find the ``appropriate linearization'' of area minimizing integral currents at flat singular points, and he developed a subsequent existence and regularity theory for these new objects;
\item[(ii)] he introduced several flexible techniques to approximate currents with Lipschitz multivalued graphs;
\item[(iii)] he developed a very intricate regularization technique to find a sufficiently smooth ``central sheet'' at possible branching singularities (the so-called ``center manifold'');
\item[(iv)] he discovered a new monotonicity formula for harmonic function (the monotonicity of the ``frequency'') which has meanwhile been used in a variety of different contexts in elliptic and parabolic partial differential equations (see e.g. \cite{GL,HS89,Log}).
\end{itemize}
Almgren's theory has been revisited by the author and Emanuele Spadaro in the series of works \cite{DS1,DS2,DS3,DS4,DS5}. Besides making the proof of Theorem \ref{t:Almgren} shorter, these works improve upon Almgren's monograph in several aspects, and moreover they have been the starting point of several further developments, which will be detailed in the next sections

Shortly after Almgren completed his 3-volumes preprint, White proved that in the case of $2$-dimensional area-minimizing currents the stratum $\mathscr{S}_0$ consists in fact of isolated points, cf. \cite{White83} (this is indeed a corollary of a more precise theorem which shows the uniqueness of tangent cones in that particular case and which will be discussed further in Section \ref{s:tangent_cones}. The program of understanding the singularities of $2$-dimensional area-minimizing currents was then completed by Chang in \cite{Chang}.

\begin{theorem}\label{t:Chang}
Let $\Sigma$ be a $2$-dimensional area minimizing current in $\mathbb R^{2+m}$. Then $\Sing_i (\Sigma)$ consists of isolated points. Moreover, for each $p\in \Sing_i (\Sigma)$ there is a neighborhood $U$ in which $\Sigma$ can be decomposed as the union of a finite number $N$ of branched minimal immersed disks $D_i$ with the following property:
\begin{itemize}
\item Each $D_i$ is an embedding except for the point $p$;
\item $D_i\cap D_j$ is either the empty set or consists only of the point $p$.
\end{itemize}
\end{theorem}

However the proof given in \cite{Chang} is strictly speaking not complete, as Chang needs the existence of a suitable generalization of Almgren's center manifold to a ``branched version''. For the latter he just gives a 4-pages sketch (cf. the appendix of \cite{Chang}), invoking suitable modifications of Almgren's statements (it must be noted that the construction of the center manifold occupies more than half of Almgren's monograph \cite{Alm00}). Based on the works \cite{DS1,DS2,DS3,DS4,DS5} the author, Spadaro, and Luca Spolaor gave a complete independent proof of the existence of a branched center manifold in \cite{DSS2}. We also developed a suitable more general counterpart of Chang's theory in the papers \cite{DSS1,DSS3,DSS4}, proving in particular the same regularity result for spherical cross-sections of area-minimizing $3$-dimensional cones and for semicalibrated $2$-dimensional currents (previous theorems in \cite{Be,BR} proved some cases of particular interest, based on the works of Rivi\`ere and Tian, see \cite{Riv,RT1,RT2}).  Almgren's dimension bound in Theorem \ref{t:Almgren} has as well been extended to semicalibrated currents by Spolaor in \cite{Spolaor}. 

\section{Interior regularity theory: minimizing currents mod $p$}

The regularity theory for area-minimizing currents mod $p$ started around the same time as the regularity theory for integral currents. As a consequence of Almgren's generalization of De Giorgi's $\varepsilon$ regularity theorem, the cases $p=2,3$ were already rather well-understood in the sixties. In particular the absence of flat singular points allowed one to infer the following Theorem (the case $p=2$ is due to Federer, cf. his pioneering work on the reduction argument \cite{Federer70}).

\begin{theorem}\label{t:mod-2}
If $\Sigma$ is an $m$-dimensional area-minimizing current mod $2$ in $\mathbb R^{m+n}$, then $\Sing_i (\Sigma) \subset \mathscr{S}_{m-2}$. If $\Sigma$ is an $m$-dimensional area-minimizing current mod $3$, then $\Sing_i (\Sigma) \subset \mathscr{S}_{m-1}$.
\end{theorem}

The case $p=2$ in codimension $1$ allows even more restrictive results (the same regularity as for integral area-minimizing currents holds and in fact locally any area-minimizing integral hypercurrent mod $2$ is the boundary of a Caccioppoli set), while in higher codimension there are indeed area-minimizing $2$-dimensional currents mod $2$ with point singularities.

For $p=3$ the union of three half-planes in $\mathbb R^3$ meeting at a common line at 120 degrees gives an obvious example for which $\mathscr{S}_{m-1}\neq \emptyset$. The beautiful result of Taylor \cite{Taylor73} gave a complete description of the interior singular set for area-minimizing $2$-dimensional currents mod $3$ in $\mathbb R^3$: locally the singular set is always diffeomorphic to the above example. The subsequent work of Simon \cite{Simon93b} on the uniqueness of cylindrical tangent cones allowed to give a suitable generalization of Taylor's result in any dimension and codimension. The final outcome is the following

\begin{theorem}\label{t:mod-3}
If $\Sigma$ is an $m$-dimensional area-minimizing current mod $3$ in $\mathbb R^{m+n}$, then $\mathscr{S}_{m-1}\setminus \mathscr{S}_{m-2}$ is an $m-1$-dimensional submanifold and at every $p\in \mathscr{S}_{m-1}\setminus \mathscr{S}_{m-2}$ there is a neighborhhod $U$ in which $\Sigma$ consists of $3$ smooth minimal surfaces meeting at $\mathscr{S}_{m-1}\cap U$ at $120$ degrees. If $\Sigma$ is in addition an hypercurrent (i.e. $n=1$), then $\mathscr{S}_{m-2}\setminus \mathscr{S}_{m-3}$ is an empty set. In particular, if $m=2$ and $n=1$, then $\Sing_i (\Sigma)$ consists of pairwise disjoint closed simple curves and pairwise disjoint simple arcs with endpoints lying in the support of the boundary of $\Sigma$.
\end{theorem}

In order to progress beyond Corollary \ref{c:Allard} for higher moduli it is necessary to either rule out flat singular points or bound their dimension. In the special case of mod $4$ hypercurrents White in \cite{White79} discovered a beautiful fact which allowed him to derive the following structural result

\begin{theorem}\label{t:mod-4}
If $\Sigma$ is an $m$-dimensional area-minimizing current mod $4$ in $\mathbb R^{m+1}$, then it can be locally decomposed, away from its boundary, in the union of two $m$-dimensional area-minimizing current mod $2$.
\end{theorem}

He also showed a converse to Theorem \ref{t:mod-4}. In particular his results imply the existence of flat singular points even for hypercurrents mod $2k$. More precisely, consider a function $u: \mathbb R^2 \supset B_1 \to \mathbb R$ which solves the minimal surface equation and, after applying a suitable translation and rotation, assume that $u(0)=0$, $\nabla u (0) =0$ and $D^2 u (0) \neq 0$. Since $\Delta u (0) =0$, it follows that the zero set of $u$ in a neighborhood of $0$ consists of 2 arcs crossing orthogonally in $0$. We can thus assume that the disk $B_r (0) \subset \mathbb R^2$ is subdivided by $\{u=0\}$ into 4 sectors $S_1, S_2, S_3, S_4$. We then consider in $\mathbb{C}_r := B_r (0)\times \mathbb R\subset \mathbb R^3$ the union of the four sectors $S_i\times \{0\}$ and of the four portions $G_i$ of the graph of $u$ lying over the respective sector $S_i$. We give to $S_i$ opposite alternating orientations and sum them to construct an integral current $S$ in $\mathbb{C}_r (0)$. Clearly $\partial S$ is formed by the four arcs which describe $\{u=0\}\times \{0\}$, suitably oriented and counted with multiplicity $2$. In particular $S$ is a cycle mod $2$. We then perform an analogous operation with the 4 portions $G_i$ of the graph of $u$ and construct a corresponding integral current $T$. By choosing the orientations correctly we can achieve that $\partial T= \partial S$. Therefore the current $\Sigma= T+S$ is a cycle mod $4$ and, according to the results in \cite{White79}, it is area-minimizing (as a cycle mod $4$). In particular $0$ is a flat singular point for $\Sigma$.

This phenomenon is typical of even moduli and indeed in his subsequent  work \cite{White84} White proved that area minimizing hypercurrents mod $2k+1$ cannot have singular flat points.

\begin{theorem}\label{t:White-odd-p}
If $\Sigma$ is an $m$-dimensional area minimizing current mod $p$ in $\mathbb R^{m+1}$ and $p$ is odd, then $\Sing_i (\Sigma) \subset \mathscr{S}_{m-1}$.
\end{theorem}

In the papers \cite{DHMS1,DHMS2} the author, Jonas Hirsch, Andrea Marchese, and Salvatore Stuvard developed a theory to bound the dimension of flat singular points of a general area-minimizing current $\Sigma$ mod $p$ (i.e. in any dimension and codimension), which implies that the Hausdorff dimension of the set of flat singular points of $\Sigma$ is at most $m-1$.

\begin{theorem}\label{t:general-p}
If $\Sigma$ is an $m$-dimensional area minimizing current mod $p$ in $\mathbb R^{m+n}$, then $\textrm{dim}_H (\Sing_i (\Sigma)) \leq m-1$.
\end{theorem}

While the latter theorem is a considerable improvement compared to what known before (aside from the cases covered by Theorems \ref{t:mod-2}, \ref{t:mod-4}, and \ref{t:White-odd-p}, in all others the best known result was that the singular set is meager, thanks to Corollary \ref{c:Allard}). Indeed the known examples would suggest that the set of flat singular points of any area minimizing current mod $p$ is at most $m-2$. The work \cite{DHMSS1} and the forthcoming one \cite{DHMSS2}, by the author, Hirsch, Marchese, Spolaor, and Stuvard give a first step towards the latter picture in codimension $1$.

\begin{theorem}\label{t:general-p-cod-1}
Let $\Sigma$ be an $m$-dimensional area minimizing current mod $p$ in $\mathbb R^{m+1}$. Then:
\begin{itemize}
\item[(a)] $\mathscr{S}_{m-1}\setminus \mathscr{S}_{m-2}$ is a $C^{1,\alpha}$ submanifold and for every $q\in \mathscr{S}_{m-1}\setminus \mathscr{S}_{m-2}$ there is a neighborhood $U$ in which $\Sigma$ consists of $p$ minimal hypersurface meeting at $\mathscr{S}_{m-1}\cap U$;
\item[(b)] At every flat singular point there is a unique tangent cone, which is a flat plane with multiplicity $\frac{p}{2}$ (in particular $p$ must be even).
\end{itemize}
\end{theorem} 
In fact, after the appearance of \cite{DHMSS1} Minder and Wickramasekera (cf. \cite{MW}) pointed out to the authors that it is possible to derive Theorem \ref{t:general-p-cod-1} directly from the theory developed in \cite{Wic}, starting from one observation in \cite{DHMSS1} concerning tangent cones in the top stratum $\mathscr{S}_{m-1}\setminus \mathscr{S}_{m-2}$ and the verification of Simon's no hole condition. In \cite{DHMSS3} Theorem \ref{t:general-p-cod-1} will be further used to confirm the conjectural picture in codimension $1$, namely to prove

\begin{theorem}\label{t:general-p-cod-1-bis}
Let $\Sigma$ be an $m$-dimensional area minimizing current mod $p$ in $\mathbb R^{m+1}$. Then $\Sing_i (\Sigma) \cap \mathscr{S}_m$ is empty for $p$ odd (as implied by Theorem \ref{t:White-odd-p}), while $\Sing_i (\Sigma)\cap \mathscr{S}_m$ has dimension at most $m-2$ for even $p$.
\end{theorem}

\section{Boundary regularity theory: minimizing integral hypercurrents}

The first boundary regularity theorem for area-minimizing integral currents $\Sigma$ was proved by Allard in his PhD thesis \cite{All-PhD} in codimension $1$. More precisely

\begin{theorem}\label{t:bdry-Allard-2}
Assume $\Sigma$ is an area-minimizing integral current of dimension $m$ in $\mathbb R^{m+1}$ and assume that 
\begin{itemize}
\item[(a)] $\partial \Sigma$ is a smooth (more precisely $C^2$) $m-1$-dimensional surface $\Gamma$ with multiplicity $1$;
\item[(b)] there is a uniformly convex smooth (more precisely $C^2$) bounded open set $U$ such that $\Gamma \subset \partial U$.
\end{itemize}
Then $\Sigma$ is smooth in a neighborhood of $\Gamma$, more precisely there is an open set $V\supset \Gamma$ such that $V\cap \Sigma$ is a smooth minimal hypersurface (with boundary) and its boundary (in $V$) is precisely $\Gamma$ (in the classical sense of differential topology).  
\end{theorem}

In fact the proof in \cite{All-PhD} contains an $\varepsilon$-regularity result which is the precursor of Theorem \ref{t:Allb}, while the Assumption (a) is combined with a suitable classification of boundary tangent cones to prove that any point $p\in \Gamma$ has density $\frac{1}{2}$. In order to remove the ``convex barrier'' of assumption (a) one needs to handle situations in which $p$ might be a ``2-sided'' boundary point. 

To illustrate the latter point, consider a $2$-dimensional plane $V$ in $\mathbb R^3$ and the two circles $\gamma_1 = \partial \mathbf{B}_1 (0) \cap V$ and $\gamma_2 = \partial \mathbf{B}_2 (0) \cap V$.  Give to $\gamma_1$ and $\gamma_2$ the ``same orientation'', so that they bound the disks $D_1 = \mathbf{B}_1 (0) \cap V$ and $D_2 = \mathbf{B}_2 (0) \cap V$, taken with the same orientation where they overlap. It can be easily shown that, if $\Gamma = \gamma_1+\gamma_2$, then $\Sigma = D_1+D_2$ is the unique area-minimizing integral current bounded by $\Gamma$. $\Sigma$ can be described as the sum of the corona $D_2\setminus D_1$, counted with multiplicity $1$, and the disk $D_1$, counted with multiplicity $2$. A point $p\in \gamma_1$ is what can be naturally called a ``2-sided'' boundary point and note that its density is $\frac{3}{2}$ (for a more rigorous definition, cf \cite{DDHM}). The regularity theory at such points is rather subtle and (in codimension $1$) it was handled in the famous work \cite{HS79} by Hardt and Simon.

\begin{theorem}\label{t:HS}
Let $\Gamma$ be a smooth oriented closed $m-1$-dimensional submanifold of $\mathbb R^{m+1}$ and let $\Sigma$ be an area-minimizing integral current whose boundary (in the sense of currents) is given by $\Gamma$ counted with multiplicity $1$. Then every point $p\in \Sigma$ is regular, namely one of the following two mutually exclusive possibilities holds:
\begin{itemize}
\item[(i)] Either the density of $\Sigma$ at $p$ is $\frac{1}{2}$ and hence the conclusion of Theorem \ref{t:Allb} applies in a neighborhood $U$ of $p$;
\item[(ii)] or the density of $\Sigma$ at $p$ is $k+\frac{1}{2}$ for some positive integer $k$; in this case there is a neighborhood $U$ of $p$ and a minimal hypersurface $\Lambda$ of $U$ without boundary such that:
\begin{itemize}
\item $\Lambda$ contains $\Gamma$;
\item $\Gamma$ subdivides $\Lambda$ in two regions $\Lambda^+$ and $\Lambda^-$;
\item $\Sigma$ in $U$ is given by $\Lambda^+$ counted with multiplicity $k+1$ and $\Lambda^-$ counted with multiplicity $k$.
\end{itemize}
\end{itemize}
\end{theorem} 

Among the many ideas introduced in \cite{HS79}, one has been highly influential in several other problems in minimal surface theory, and it is the so called Hardt-Simon inequality. In a nutshell the Hardt-Simon inequality makes a clever use of the remainder in the monotonicity formula (namely the precise expression for the quantity $\frac{d}{dr} \frac{\textrm{Vol}^m (\Sigma \cap \mathbf{B}_r)}{r^m}$) in order to infer nontrivial information on the graphical approximation of $\Sigma$ at small scales. 

While we have stated Theorems \ref{t:bdry-Allard-2} and \ref{t:HS} as ``global theorems'', suitable local versions of them are also valid, and in fact the very nature of the main arguments is completely local.

\section{Boundary regularity theory: minimizing integral currents with smooth boundaries of multiplicity $1$}

In his fundamental boundary regularity paper \cite{Allard75} Allard noticed that Theorem \ref{t:Allb} can be used to generalize the conclusion of Theorem \ref{t:bdry-Allard-2} to all codimensions.

\begin{theorem}\label{t:bdry-Allard-3}
Let $\Gamma$ be a smooth $m-1$-dimensional closed oriented submanifold of $\mathbb R^{m+n}$ and let $U$ be a bounded smooth uniformly convex set such that $\Gamma \subset \partial U$. 
Then any area-minimizing integral current $\Sigma$ whose boundary is given by $\Gamma$ (counted with multiplicity $1$) is smooth in a neighborhood of $\Gamma$, in the sense of the conclusion of Theorem \ref{t:bdry-Allard-2}.
\end{theorem}

Again, a local version of the above theorem holds as well and in fact,in order to conclude that a boundary point $p$ is regular and one-sided in the sense of Theorem \ref{t:HS}(i), it suffices to find a uniformly convex ``barrier'' which touches $\Gamma$ at $p$ and so that $\Sigma$ lies (locally) on one side of it, cf. \cite{Hardt}. A simple argument furnishes such a barrier for any smooth $\Gamma\subset \mathbb R^{m+n}$: for instance one could consider the smallest closed ball containing $\Gamma$. It then follows that under the mere assumption that $\Gamma$ is sufficiently smooth, an area-minimizing current bounding $\Gamma$ (taken with multiplicity $1$) has always at least one boundary regular point. 

Up until recently nothing more was known, except that in codimension higher than $1$ singular boundary points are certainly possible. A simple example is given by the union of a smooth simple curve $\gamma_1\subset \{x_1=x_2=0\}\subset \mathbb R^4$ containing the origin and a smooth simple curve $\gamma_2\subset \{x_3=x_4=0\}\subset \mathbb R^4$ which {\em does not} contain the origin. This union bounds an area-minimizing $2$-dimensional integral current for which $0$ is a boundary singular point. Moreover, since in a general Riemannian manifold the barrier argument outlined in the previous paragraph is not available, even in the simplest case of a smooth simple closed curve in a closed smooth Riemannian $4$-manifold $M$, the results outlined so far could not exclude the possibility that {\em all} boundary points of an area-minimizing current $\Sigma \subset M$ bounding $\Gamma$ are singular. 

As in the case of Theorem \ref{t:HS} the main difficulty in removing the convex barrier assumption is the possibility that boundary points have density larger than $\frac{1}{2}$. And as in the case of Theorem \ref{t:Almgren} the most problematic issue is that unfortunately the existence of a flat tangent cone at the boundary does not guarantee regularity: flat boundary singular points exist as soon as the codimension is larger than $1$, cf. \cite{DDHM}. In \cite{DDHM} the author, Guido De Philippis, Hirsch, and Annalisa Massaccesi were able to develop a suitable ``Almgren-type'' regularity theory for boundary points, building on a previous important step of Hirsch \cite{Hirsch}. In particular we proved the following

\begin{theorem}\label{t:DDHM}
Let $\Gamma$ be a smooth closed oriented $m-1$-dimensional submanifold of $\mathbb R^{m+n}$ and let $\Sigma$ be an area minimizing integral current whose boundary is given by $\Gamma$ taken with multiplicity $1$. Then the set of boundary regular points, understood as points where one of the two alternatives (i) and (ii) of Theorem \ref{t:HS} hold, is a dense relatively open subset of $\Gamma$.
\end{theorem} 

While Theorem \ref{t:DDHM} might look very far from optimal, it turns out that a naive counterpart of the bound of the dimension of the interior singular set is in fact false. In \cite{DDHM} we prove also the following.

\begin{theorem}\label{t:DDHM-2}
There is a smooth $1$-dimensional embedded submanifold of $\mathbb R^4$ which bounds an area-minimizing current $\Sigma$ of $\mathbb R^4$ whose boundary singular set has Hausdorff dimension $1$.
\end{theorem}

Theorem \ref{t:DDHM-2} leaves open the possibility that at least the set of boundary singular points has zero $m-1$-dimensional Hausdorff measure and that it has dimension $m-2$ if the boundary is real analytic. We also caution the reader that a less restrictive definition of boundary regular point might restrict the size of boundary singularities even in the $C^\infty$ case. For a more detailed discussion of all these possibilities we refer the reader to Section \ref{s:open-problems}. However, that boundary regularity is subtle is also witnessed by the following example of the author, De Philippis, and Hirsch (cf. \cite{DDH}).

\begin{theorem}\label{t:DDH}
There is a smooth closed $4$-dimensional Riemannian manifold $M$ and a smooth simple closed curve $\Gamma\subset M$ which bounds a unique area-minimizing $2$-dimensional current $\Sigma$ which is smooth in $M\setminus \Gamma$ and whose first homology group is infinite-dimensional. In fact $\Sigma$ is smooth except at a single point $p\in \Gamma$. 
\end{theorem}

\section{Boundary regularity theory: minimizing integral currents with smooth boundaries of higher multiplicity}

In the previous sections we examined the boundary regularity of area-minimizing integral currents under the assumption that the multiplicity of the boundary is $1$. A rather intriguing and widely open question, already raised by Allard in his PhD thesis \cite{All-PhD}, is what happens when the multiplicity is an integer larger than $1$ (the fact that it {\em must} be an integer is of course a consequence of the integrality assumption, but we also remind the reader that when $T$ is integer rectifiable and $\partial T$ has finite mass, $\partial T$ is necessarily integer rectifiable, cf. \cite{Federer}). 

The problem raised by Allard in \cite{All-PhD} is highlighted again by White in \cite{Collection}. In the same reference White observes also that, thanks to the decomposition theorem for area-minimizing hypercurrents, if $\Sigma$ is an $m$-dimensional area-minimizing integral current in $\mathbb R^{m+1}$ whose boundary is a smooth submanifold $\Gamma$ counted with multiplicity $Q>1$, then $\Sigma$ can be decomposed into the sum of $Q$ area-minimizing integral currents whose boundary is $\Gamma$ counted with multiplicity $1$ also . In particular we are in the position of applying the Hardt-Simon Theorem \ref{t:HS} to each element of the decomposition. While this is the same ``codimension 1 phenomenon'' that rules out flat singular points in the interior for area-minimizing integral hypercurrents, we pause a moment to make one important remark. It is well known that there are smooth $m-1$-dimensional oriented closed submanifolds of $\mathbb R^{m+1}$ that bound more than one area-minimizing integral current. This is already the case for smooth simple closed curves $\Gamma$  in $\partial \mathbf{B}_1 \subset \mathbb R^3$. Consider in particular one such $\Gamma$ and let $\Sigma_1$ and $\Sigma_2$ be two area-minimizing integral $2$-dimensional currents which bound $\Gamma$ (with multiplicity $1$). Thanks to the above decomposition theorem $\Sigma = \Sigma_1+\Sigma_2$ is an area-minimizing current which bounds a double copy of $\Gamma$\footnote{It is one of the most beautiful discoveries of geometric measure theory that this conclusion is in general false in higher codimension. In particular, following the pioneering work of L. C. Young \cite{Young63}, there are several constructions of smooth simple curves $\Gamma$ in $\mathbb R^4$ with the following remarkable property. If we let $m (\Gamma)$ be the mass of an area-minimizing $2$-dimensional current which bounds one copy of $\Gamma$ and $m (2\Gamma)$ the mass of an area-minimizing integral current which bounds two copies of $\Gamma$, then $m (2\Gamma) < 2 m (\Gamma)$.}. By the interior regularity theory $\Sigma_1$ and $\Sigma_2$ have no interior point in common. Therefore, by the Hopf boundary lemma, there is no boundary point at $\Gamma$ in which $\Sigma_1$ and $\Sigma_2$ have the same tangent: $\Sigma_1$ and $\Sigma_2$ meet at every point of $\Gamma$ {\em transversally}. 

In light of the above example,  it seems sensible to give the following definition of boundary regular point.

\begin{definition}\label{d:bdry-regular}
Assume that $\Gamma\subset \mathbb R^{m+n}$ is a smooth oriented $m-1$-dimensional submanifold and that $Q$ is a positive integer. Let $\Sigma$ be an area-minimizing integral current in $\mathbb R^{m+n}$ whose boundary is given by $Q$ copies of $\Gamma$. $p\in \Gamma$ is a regular boundary point if one of the following two alternatives occur in some neighborhood $U$ of $p$:
\begin{itemize}
\item[(i)] There are $N$ positive integers $k_i$ with $\sum_i k_i =Q$ and $N$ smooth minimal surfaces $\Lambda_i$ in $U$ with boundary $\Gamma$ such that $\Sigma \cap U = \sum_i k_i \Lambda_i$, and each distinct pair $\Lambda_i$ and $\Lambda_j$ meet transversally at $p$; 
\item[(ii)] There is a minimal surface $\Lambda$ in $U$ without boundary, which contains $\Gamma$: the latter subdivides $\Lambda$ in two regions $\Lambda^+$ and $\Lambda^-$ and $\Sigma \cap U = (Q+k) \Lambda^+ + Q\Lambda^-$ for some positive integer $k$.
\end{itemize}
\end{definition} 

In particular the discussion above reduces the following statement to a mere corollary of Theorem \ref{t:HS}:

\begin{corollary}\label{c:HS2}
Let $\Gamma$ be a smooth oriented closed $m-1$-dimensional submanifold of $\mathbb R^{m+1}$, $Q$ be a positive integer, and $\Sigma$ an area-minimizing integral current with $\partial \Sigma = Q \Gamma$. Then every boundary point $p\in \Gamma$ is regular in the sense of Definition \ref{d:bdry-regular}.  
\end{corollary}

The boundary regularity theory for $Q>1$ and in codimension larger than $1$ is widely open. A very first preliminary result, which is a counterpart of Theorem \ref{t:bdry-Allard-2} for $2$-dimensional area-minimizing currents, has been proved very recently by the author, Stefano Nardulli, and Simone Steinbr\"uchel in \cite{DNS1,DNS2}, building in part upon the theory developed in \cite{DDHM} and the paper \cite{HM}. 

\begin{theorem}\label{t:DNS}
Consider a smooth $1$-dimensional closed submanifold $\Gamma$ of $\mathbb R^{2+n}$ and assume that there is a bounded smooth uniformly convex open set $U$ such that $\Gamma \subset \partial U$. Let $\Sigma$ be an area-minimizing $2$-dimensional integral current such that $\partial \Sigma = Q \Gamma$ for some integer $Q$. Then every point $p\in \Gamma$ is a boundary regular point and moreover alternative (i) in Definition \ref{d:bdry-regular} holds at every such point. 
\end{theorem}

\section{Uniqueness of tangent cones}\label{s:tangent_cones}

One major open question in the regularity theory of minimal submanifolds, which has attracted the attention of a large number of researchers since the dawn of geometric measure theory, is the uniqueness of tangent cones. This amounts to the question of whether there is at every point $p$ a unique limit for the rescalings $\frac{\Sigma-p}{r}$ of the minimal submanifold $\Sigma$. In some situations the question is intimately connected to the understanding of the regularity properties of the various strata $\mathscr{S}_k\setminus \mathscr{S}_{k-1}$. For instance the pioneering works of Taylor \cite{Taylor73,Taylor76} leading to the Theorems  \ref{t:Taylor} and \ref{t:mod-3} can be reduced to suitable uniqueness statements for the relevant tangent cones. 

The most striking result in the area is the celebrated theorem of Simon in \cite{Simon83}.

\begin{theorem}\label{t:Simon}
Let $\Sigma$ be a stationary integral varifold and assume that the spherical cross section of one tangent cone $\Sigma_0$ at an interior point $p$ of $\Sigma$ is a regular submanifold of $\partial \mathbf{B}_1$ with multiplicity $1$. Then $\Sigma_0$ is the unique tangent cone to $\Sigma$ at $p$. 
\end{theorem}

Once again a similar uniqueness theorem is widely open when the multiplicity of the cross section is allowed to take multiplicity higher than $1$, except for some lucky situations in which the case of higher multiplicity can be reduced to that of multiplicity $1$. Two notable examples are that of area-minimizing hypercurrents and that of area-minimizing currents mod $2$

\begin{corollary}\label{c:Loj-Simon}
Let $\Sigma$ be an area-minimizing $m$-dimensional integral current in $\mathbb R^{m+1}$ or an area-minimizing $m$-dimensional current mod $2$ in $\mathbb R^{m+n}$ and assume that $p$ is an interior point at which one tangent cone $\Sigma_0$ has smooth spherical cross section. Then $\Sigma_0$ is the unique tangent cone at $p$. 
\end{corollary}

Even in the case of multiplicity $1$, the uniqueness of tangent cones whose spherical cross section is not smooth is a much more subtle issue. Before discussing it, we wish to introduce a suitable concept which has played a pivotal role in many contexts. Let $\Sigma_0$ be a stationary varifold which is a cone with smooth cross section $\Gamma_0$, taken with multiplicity $1$. It is then well-known that $\Gamma_0$ is a minimal submanifold of the sphere $\partial \mathbf{B}_1$. If $\Sigma_k$ is a sequence of cones converging to $\Sigma_0$, with cross section $\Gamma_k$, up to extraction of a subsequence $\Gamma_k$ is the graph of a solution of a suitable linear elliptic PDE over $\Gamma_0$ plus higher order terms. Such solutions are called Jacobi fields in the literature, and they are the higher-dimensional counterpart of the classical Jacobi fields on geodesics. A Jacobi field $u$ is called ``integrable'' if there is a sequence $\Gamma_k$ of minimal submanifolds of $\partial \mathbf{B}_1$ covering to $\Gamma_0$ which generates $u$ as outlined above. Prior to Theorem \ref{t:Simon} Allard and Almgren in \cite{All-Alm} proved the following important result

\begin{theorem}\label{t:Allard-Almgren}
Let $\Sigma$ be a stationary integral varifold, let $p$ be an interior point and assume that a tangent cone $\Sigma_0$ to $\Sigma$ at $p$ satisfies the following two properties:
\begin{itemize}
\item[(i)] The spherical cross section $\Gamma_0$ of $\Sigma_0$ is smooth and is taken with multiplicity $1$;
\item[(ii)] Every Jacobi field of $\Gamma_0$ is integrable.
\end{itemize}
Then $\Sigma_0$ is the unique tangent cone to $\Sigma$ at $p$ and moreover the rescalings $\frac{\Sigma -p}{r}$ converge to $\Sigma_0$ with a polynomial rate. 
\end{theorem}

The integrability condition (ii) of Theorem \ref{t:Allard-Almgren} has several drawbacks. In order to verify it one must know rather explicitly the cross section $\Gamma_0$. But even in the cases in which $\Gamma_0$ is known and has a rather simple formula, verifying the condition 
is in general quite hard (in particular it requires a classification result for all the solutions of some particular elliptic PDE). Last but not least, there are examples in which it does not hold, see \cite{AS}, and in which the convergence rate of $\frac{\Sigma-p}{r}$ to $\Sigma_0$ is just logarithmic. The powerful approach of Simon to Theorem \ref{t:Simon} avoids any discussion of the integrability of the Jacobi vector fields thanks to his realization that the convergence of $\frac{\Sigma-p}{r}$ to $\Sigma_0$ can be reduced to an infinite-dimensional version of a classical result of Lojasiewicz for finite-dimensional gradient flows. The corresponding ``Lojasiewicz-Simon inequality'' has been widely used to study the convergence of parabolic PDEs to a unique steady state and the uniqueness of model singularities in other geometric variational problems. 

Coming back to cones whose spherical cross-sections are not smooth, a particularly simple subclass are called ``cylindrical tangent cones''. In his notable investigation \cite{Simon93b} Simon has been able to prove a useful generalization of the Allard-Almgren Theorem \ref{t:Allard-Almgren}.

\begin{theorem}\label{t:Simon-cylindrical}
Let $\Sigma$ be a stationary integral varifold, let $p$ be an interior point and assume that a tangent cone $\Sigma_0$ to $\Sigma$ at $p$ satisfies the following structural properties:
\begin{itemize}
\item[(i)] $\Sigma_0 = V \times \Lambda_0$ for some minimal cone $\Lambda_0$ and some linear subspace $V$;
\item[(ii)] The spherical cross section $\Gamma_0$ of $\Lambda_0$ is smooth and is taken with multiplicity $1$;
\item[(iii)] Every Jacobi field of $\Gamma_0$ is integrable;
\item[(iv)] The following ``no hole condition'' holds for a sufficiently small $\delta (\Gamma_0)$: provided $\frac{\Sigma-p}{r}$ is sufficiently close to $\Sigma_0$, every $\mathbf{B}_\delta (q)$ with $q\in \mathbf{B}_1 \cap V$ contains a point $x$ of density $\Theta (\Sigma, x)$ larger than $\Theta (\Sigma, p)-\delta$.
\end{itemize}
Then $\Sigma_0$ is the unique tangent cone to $\Sigma$ at $p$ and moreover the rescalings $\frac{\Sigma -p}{r}$ converge to $\Sigma_0$ with a polynomial rate. 
\end{theorem} 

Quite a few of the structural results for singular strata mentioned in the previous sections depend heavily on the above result (or can be deduced from it). A notable exception is the uniqueness theorem of Taylor which underlines the second conclusion of Theorem \ref{t:Taylor} (and the higher-dimensional counterpart in \cite{COS}).  The latter is in fact derived through a direct epiperimetric inequality \`a la Reifenberg. 

One major drawback of the approaches to Theorems \ref{t:Simon}, \ref{t:Allard-Almgren}, and \ref{t:Simon-cylindrical} is that the underlying PDE arguments relie on heavily on the $\varepsilon$-regularity result of Allard, namely Theorem \ref{t:Allard} (or on some other analogous results). For instance in Theorem \ref{t:Simon} the multiplicity $1$ assumption and the regularity of $\Gamma_0$ allow to conclude that $\partial \mathbf{B}_1 \cap \frac{\Sigma-p}{r}$ is a smooth graph over $\Gamma_0$. On the contrary, the epiperimetric inequality, which is based on exhibiting a suitable competitor, can be applied in situations where the cross section is irregular or taken with higher multiplicity. On the other hand its applicability is limited to minimizers. Up until recently, another obvious objection to a wider plausibility of an epiperimetric inequality \`a la Reifenberg is that it immediately implies a polynomial decay rate, which is known to be false in general cf. \cite{AS}. However the recent paper of Colombo, Spolaor, and Velichkov \cite{CSV} shows that Theorem \ref{t:Simon} can be recovered (and in fact generalized to suitable ``quasi-minima'') through a suitable generalization of Reifenberg's epiperimetric (it must be noted that the proof of the latter is nonetheless achieved using the Lojasiewicz-Simon inequality).

An important case in which an epiperimetric inequality can be proved and used effectively to prove uniqueness of tangent cones (while a ``PDE-approach'' has not yet been given) is that of $2$-dimensional area-minimizing currents at interior points. In particular White in \cite{White83} proved 

\begin{theorem}\label{t:White}
Let $\Sigma$ be a $2$-dimensional area minimizing integral current in $\mathbb R^{2+n}$. Then the tangent cone to $\Sigma$ is unique at every interior point $p$. 
\end{theorem}

A counterpart to Theorem \ref{t:White} has been shown by Hirsch and Marini in \cite{HM} at smooth boundaries taken with multiplicity $1$. However, as noticed in \cite{DNS1}, the proof of Hirsch and Marini can easily be adapated to the case of smooth boundaries with arbitrary multiplicities, thus giving a complete result for $2$-dimensional area-minimizing integral currents in any dimension and codimension. 

\section{Open problems}\label{s:open-problems}

In this section we will collect some open questions. I wish to emphasize that the selection given here by no means exhaust the interesting open problems in the area, but it rather reflects a personal choice of the author.

\subsection{Stationary and stable varifolds} Perhaps the most intriguing question is whether it is possible to improve Corollary \ref{c:Allard} in any situation which is not the trivial one of $1$-dimensional stationary varifolds. The most modest goal would be to show that the singular set of stationary $2$-dimensional integral varifolds in $\mathbb R^3$ has zero $2$-dimensional Hausdorff measure. In general there is no example of singular stationary $m$-dimensional varifolds (in any codimension) for which the singular set has dimension larger than  $m-1$. 

In the case of stable varifolds of codimension $1$ the deep theory of Wickramasekera developed in \cite{Wic} (see also \cite{Minter} for some further progress) makes one hope that in the future some final unconditional structural result might be at hand. A coronation of the efforts in the area would be a theorem which proves that $\mathscr{S}_{m-1}\setminus \mathscr{S}_{m-2}$ is a $C^{1,\alpha}$ $m-1$-dimensional submanifold, while the set of flat singularities is $m-2$-rectifiable. The latter statement seems to be reachable in the very particular case of area-minimizing hypercurrents mod $p$. 

A widely open problem is whether stability allows one to go beyond Allard's conclusion in codimension higher than $1$. It is quite baffling that no further regularity information has been concluded thus far for stable varifolds as soon as the codimension is larger than $1$. A particularly intriguing case would be that of $2$-dimensional stable varifolds, already in $\mathbb R^4$. A first question in that direction is whether some counterpart of the Schoen-Simon-Yau estimates and hence a corresponding compactness theorem holds for classical (possibly branched) minimal $2$-dimensional surfaces in, say, $\mathbb R^4$. In other words, assume that $\Sigma_k$ is such a sequence in some bounded open set $U\subset \mathbb R^4$, that the area of $\Sigma_k$ is uniformly bounded, and that each $\Sigma_k$ is stable. Is it possible to extract a subsequence which is converging (in the varifold sense) to a classical stable (possibly branched) immersed minimal surface? What if we restrict further $\Sigma_k$ and ask that they are embedded except for a finite number of branching singularities? Note that Theorem \ref{t:Chang} does imply the desired conclusion if each $\Sigma_k$ can be oriented so to give an area-minimizing integral $2$-dimensional current.

\subsection{Singularities of area-minimizing integral hypercurrents} Area-minimizing integral currents of dimension $m$ in $\mathbb R^{m+1}$ are the objects for which we have the strongest regularity theory. Is it possible to prove more about the structure of the singular set? In particular, is it true that $\mathscr{S}_{m-7}\setminus \mathscr{S}_{m-8}$ is a $C^{1,\alpha}$ submanifold, or rather are there examples (as the recent stable minimal hypersurfaces in some Riemannian manifolds given by Simon in \cite{Simon21b}) for which $\mathscr{S}_{m-7}$ has a fractal Hausdorff dimension $m-8 < \alpha < m-7$? Does it make a difference if the ambient is a smooth Riemannian manifold rather than Euclidean space?

A closely related question is whether the no-hole condition of (iv) in Theorem \ref{t:Simon-cylindrical} can be violated at some point $p$ of an area-minimizing hypercurrent (in the Euclidean space or in a general smooth Riemannian ambient). This is indeed the case for some points in the examples of stable minimal hypersurfaces constructed in \cite{Simon21b}, while in the Euclidean space a completely different example has been given by Gabor Sz\'ekelyhidi in \cite{Sz2}.

\subsection{Singularities of area-minimizing integral currents in codimension higher than $1$} It is very tempting to conjecture that for $m\geq 3$ Almgren's partial regularity theorem can improved to say that the singular set of any area-minimizing integral $m$-dimensional current in $\mathbb R^{m+n}$ is $m-2$-rectifiable. This problem seems intimately linked to the ``simplest'' open case of uniqueness of tangent cones for area-minimizing currents in codimension $n\geq 2$:
\begin{itemize}
\item Consider an area-minimizing integral current $\Sigma$ of dimension $m$ in $\mathbb R^{m+n}$ and let $p\in \Sing_i (\Sigma)$ be a point where one tangent cone is flat. Is the latter the {\em unique} tangent cone to $\Sigma$ at $p$?
\end{itemize} 
The forthcoming work \cite{DS} seems to suggest that a positive answer to the latter question, together with the additional information that the convergence rate is polynomial, would imply $m-2$ rectifiability of $\Sing_i (\Sigma)$.

On the other hand the works \cite{L1,L2} suggest that further structural results cannot be expected, at least not in general smooth ambient manifolds, and instead there are $3$-dimensional area-minimizing integral currents in closed smooth Riemannian manifolds whose singular sets have any preassigned Hausdorff dimension $\alpha \in (0,1)$.

\subsection{Singularities of area-minimizing currents mod $p$} As already mentioned, in the works \cite{DHMSS1,DHMSS2,DHMSS3} (see also \cite{MW,Wic}) we plan to show that, for an $m$-dimensional area-minimizing current mod $p$ in $\mathbb R^{m+1}$, the stratum $\mathscr{S}^{m-1}\setminus \mathscr{S}^{m-2}$ is a $C^{1,\alpha}$ $m-1$-dimensional submanifold, while the set of flat singular points has dimension at most $m-2$. In fact it is expected that the latter is $m-2$-rectifiable. The same properties could be expected in higher codimension, but the problem poses considerable difficulties. Moreover, the author does not know examples in which the stratum $\mathscr{S}^{m-2}\setminus \mathscr{S}^{m-3}$ is nonempty. In that respect the most basic question is whether there is any counterpart of Taylor's theorem for the case $p=3$: is there any $2$-dimensional area-minimizing cone mod $p$ in $\mathbb R^3$ which is not invariant under some translation? \cite{Taylor73} and \cite{White79} imply that the answer is no for $p=3$ and $4$ (while it is a simple exercise to see that it is no for $p=2$ as well, since it reduces to the case of integral currents). 

\subsection{Boundary regularity of area-minimizing integral currents at multiplicity $1$ boundaries} Is it possible to improve Theorem \ref{t:DDHM} and show that for general smooth $\Gamma$ the set of boundary singular points has zero Hausdorff $m-1$-dimensional measure? It must also be noted that, in the examples of Theorem \ref{t:DDHM-2} given by the argument of \cite{DDHM}, most of the boundary singular points $p$'s are of ``crossing types'', i.e. in some neighborhood $U$ of such $p$'s the area-minimizing current can be decomposed in one area-minimizing current which takes the boundary $\Gamma$ smoothly and a second one which is area-minimizing and has no boundary (but includes $p$ in its support). In particular the following two conjectures seem likely:
\begin{itemize}
\item Boundary singularities of non-crossing type have a much lower dimension (according to the examples the best we can hope is $m-2$).
\item Since crossing type singularities have necessarily dimension $m-2$ when $\Gamma$ (and the ambient Riemannian manifold) is real analytic, the whole boundary singular set has dimension at most $m-2$ under the latter assumption.
\end{itemize} 
In fact the following elegant conjecture is due to White in \cite{White96}.

\begin{conjecture}
Let $\Gamma\subset \mathbb R^{2+n}$ be a simple closed real-analytic curve and $\Sigma$ an area-minimizing integral current such that $\partial \Sigma = \Gamma$. Then the union of the boundary and interior singular points of $\Sigma$ is discrete. In particular:
\begin{itemize}
\item The ``overall singular set'' is finite, 
\item $\Sigma$ has finite genus $g$
\item and it is a classical Douglas-Rado solution of the Plateau problem among surfaces of genus $g$. 
\end{itemize}
\end{conjecture}

\subsection{Boundary regularity of area-minimizing integral currents at boundaries with higher multiplicity} It is tempting to conjecture that Theorem \ref{t:DNS} holds for $m$-dimensional integral currents for $m\geq 2$, but in reality the situation might be more complicated. Otherwise a more modest expectation is that for general $m$, under the assumptions of Theorem \ref{t:DNS}, the boundary singular set has dimension at most $m-3$. Nothing is known in the case of a general integral multiplicity $Q$ and a general boundary $\Gamma$, i.e. without the assumption that there is a ``convex barrier'' at (a portion of) $\Gamma$. One might expect that the counterpart of Theorem \ref{t:DDHM} holds for general multiplicities $Q\geq 1$.

\subsection{Uniqueness of tangent cones} The uniqueness of interior tangent cones when the multiplicity of the cross section is larger than $1$ is widely open. As already mentioned, the most striking case is that of flat singular points, i.e. points at which at least one tangent cone is a plane with higher multiplicity, but the generalized minimal surface is not regular. This problem is open for integral area-minimizing currents of dimension $m\geq 3$ in codimension larger than $n\geq 2$, but it is also open for stationary and stable varifolds in dimension $m\geq 2$ and codimension $1$.

It is also widely open whether Simon's Theorem \ref{t:Simon-cylindrical} can be improved. In particular, can one drop the ``no-hole condition'' (iv) or the integrability condition (iii), at least for some suitable subclass of stationary varifolds? Some situations in which the ``no-hole condition'' can be dropped are given in \cite{Simon94}, while the recent work \cite{Sz1} is the first, to the author knowledge, in which the uniqueness of the cylindrical cone is proved for one example in which both conditions (iii) and (iv) in Theorem \ref{t:Simon-cylindrical} can be dropped. 

We finish this survey by mentioning that the following very innocent question is still open (even in the case $m=3$ and $n+2=5$):
\begin{itemize}
\item Consider an $m$-dimensional area-minimizing integral current $\Sigma$ in $\mathbb R^{m+n}$ with $m\geq 3$ and $n\geq 2$. Assume that one tangent cone $\Sigma_0$ at some point $p\in \Sing_i (\Sigma)$ is the union of two distinct linear planes counted both with multiplicity $1$. Is $\Sigma_0$ the unique tangent cone to $\Sigma$ at $p$?
\end{itemize}

\end{document}